\documentclass[12pt]{article}
\usepackage{full page,
latexsym, amssymb, amscd, amsthm, graphicx, amsmath, 
}

\begin{document} 

\title{\bf Logarithmic Intertwining Operators and Genus-One Correlation Functions} 
\author{Francesco Fiordalisi} 
\date{}

\bibliographystyle{alpha}
\maketitle
 
\newtheorem{thm}{Theorem}[section]
\newtheorem{defn}[thm]{Definition}
\newtheorem{prop}[thm]{Proposition}
\newtheorem{cor}[thm]{Corollary}
\newtheorem{lemma}[thm]{Lemma}
\newtheorem{rema}[thm]{Remark}
\newtheorem{app}[thm]{Application}
\newtheorem{prob}[thm]{Problem}
\newtheorem{conv}[thm]{Convention}
\newtheorem{conj}[thm]{Conjecture}
\newtheorem{cond}[thm]{Condition}
\newtheorem{exam}[thm]{Example}
\newtheorem{assum}[thm]{Assumption}
\newtheorem{nota}[thm]{Notation}
\newcommand{\halmos}{\rule{1ex}{1.4ex}}
\newcommand{\pfbox}{\hspace*{\fill}\mbox{$\halmos$}}

\newcommand{\nn}{\nonumber \\}
\newcommand{\bn}{\displaybreak[1] \nn}

\newcommand{\rad}{\mbox{\rm rad}}
\newcommand{\res}{\mbox{\rm Res}}
\newcommand{\ord}{\mbox{\rm ord}}
\renewcommand{\hom}{\mbox{\rm Hom}}
\newcommand{\edo}{\mbox{\rm End}\ }
\newcommand{\pf}{{\it Proof.}\hspace{2ex}}
\newcommand{\epf}{\hspace*{\fill}\mbox{$\halmos$}}
\newcommand{\epfv}{\hspace*{\fill}\mbox{$\halmos$}\vspace{1em}}
\newcommand{\epfe}{\hspace{2em}\halmos}
\newcommand{\nord}{\genfrac{}{}{0pt}{2}{\circ}{\circ}}
\newcommand{\wt}{\mbox{\rm wt}\ }
\newcommand{\swt}{\mbox{\rm {\scriptsize wt}}\ }
\newcommand{\lwt}{\mbox{\rm wt}^{L}\;}
\newcommand{\rwt}{\mbox{\rm wt}^{R}\;}
\newcommand{\slwt}{\mbox{\rm {\scriptsize wt}}^{L}\,}
\newcommand{\srwt}{\mbox{\rm {\scriptsize wt}}^{R}\,}
\newcommand{\clr}{\mbox{\rm clr}\ }
\newcommand{\tr}{\mbox{\rm tr}}
\newcommand{\Tr}{\mbox{\rm tr}^\phi_{\tilde W_n}}
\newcommand{\C}{\mathbb{C}}
\newcommand{\Z}{\mathbb{Z}}
\newcommand{\R}{\mathbb{R}}
\newcommand{\Q}{\mathbb{Q}}
\newcommand{\N}{\mathbb{N}}
\newcommand{\CN}{\mathcal{N}}
\newcommand{\D}{\mathcal{D}}
\newcommand{\F}{\mathcal{F}}
\newcommand{\I}{\mathcal{I}}
\newcommand{\sO}{\mathcal{O}}
\newcommand{\sQ}{\mathcal{Q}}
\newcommand{\V}{\mathcal{V}}
\newcommand{\U}{\mathcal{U}}
\newcommand{\one}{\mathbf{1}}
\newcommand{\ce}{\mathbf{k}}
\newcommand{\hz}{\hat{\mathfrak{h}}_{\Z +\frac{1}{2}}}
\newcommand{\hzt}{\hat{\mathfrak{h}}_{\Z}}
\newcommand{\BY}{\mathbb{Y}}
\newcommand{\ds}{\displaystyle}
\newcommand{\sfrac}[2]{\genfrac{}{}{1pt}{1}{#1}{#2}}
\newcommand{\ssfrac}[2]{\genfrac{}{}{1pt}{2}{#1}{#2}}

        \newcommand{\ba}{\begin{array}}
        \newcommand{\ea}{\end{array}}
        \newcommand{\be}{\begin{equation}}
        \newcommand{\ee}{\end{equation}}
        \newcommand{\bea}{\begin{eqnarray}}
        \newcommand{\eea}{\end{eqnarray}}
        \newcommand{\lbar}{\bigg\vert}
        \newcommand{\p}{\partial}
        \newcommand{\dps}{\displaystyle}
        \newcommand{\bra}{\langle}
        \newcommand{\ket}{\rangle}

        \newcommand{\ob}{{\rm ob}\,}
        \renewcommand{\hom}{{\rm Hom}}

\newcommand{\A}{\mathcal{A}}
\newcommand{\Y}{\mathcal{Y}}
\newcommand{\gY}[1]{\Y_{#1}(\U(x_{#1})w_{#1},x_{#1})}
\newcommand{\gYq}[1]{\Y_{#1}(\U(q_{z_{#1}})w_{#1}, q_{z_{#1}})}
\newcommand{\qL}{q^{L(0)}}
\newcommand{\End}{\mathrm{End}}

\newcommand{\der}[2]{\frac{\partial^{#1}}{\partial #2^{#1}}}

\newcommand{\dlt}[3]{#1 ^{-1}\delta \bigg( \frac{#2 #3 }{#1 }\bigg) }

\newcommand{\dlti}[3]{#1 \delta \bigg( \frac{#2 #3 }{#1 ^{-1}}\bigg) }

\renewcommand{\theequation}{\thesection.\arabic{equation}}
\renewcommand{\thethm}{\thesection.\arabic{thm}}
\setcounter{equation}{0} \setcounter{thm}{0}
\setcounter{section}{-1}

 \makeatletter
\newlength{\@pxlwd} \newlength{\@rulewd} \newlength{\@pxlht}
\catcode`.=\active \catcode`B=\active \catcode`:=\active \catcode`|=\active
\def\sprite#1(#2,#3)[#4,#5]{
   \edef\@sprbox{\expandafter\@cdr\string#1\@nil @box}
   \expandafter\newsavebox\csname\@sprbox\endcsname
   \edef#1{\expandafter\usebox\csname\@sprbox\endcsname}
   \expandafter\setbox\csname\@sprbox\endcsname =\hbox\bgroup
   \vbox\bgroup
  \catcode`.=\active\catcode`B=\active\catcode`:=\active\catcode`|=\active
      \@pxlwd=#4 \divide\@pxlwd by #3 \@rulewd=\@pxlwd
      \@pxlht=#5 \divide\@pxlht by #2
      \def .{\hskip \@pxlwd \ignorespaces}
      \def B{\@ifnextchar B{\advance\@rulewd by \@pxlwd}{\vrule
         height \@pxlht width \@rulewd depth 0 pt \@rulewd=\@pxlwd}}
      \def :{\hbox\bgroup\vrule height \@pxlht width 0pt depth
0pt\ignorespaces}
      \def |{\vrule height \@pxlht width 0pt depth 0pt\egroup
         \prevdepth= -1000 pt}
   }
\def\endsprite{\egroup\egroup}
\catcode`.=12 \catcode`B=11 \catcode`:=12 \catcode`|=12\relax
\makeatother

\def\hboxtr{\FormOfHboxtr} 
\sprite{\FormOfHboxtr}(25,25)[0.5 em, 1.2 ex] 

:BBBBBBBBBBBBBBBBBBBBBBBBB |
:BB......................B |
:B.B.....................B |
:B..B....................B |
:B...B...................B |
:B....B..................B |
:B.....B.................B |
:B......B................B |
:B.......B...............B |
:B........B..............B |
:B.........B.............B |
:B..........B............B |
:B...........B...........B |
:B............B..........B |
:B.............B.........B |
:B..............B........B |
:B...............B.......B |
:B................B......B |
:B.................B.....B |
:B..................B....B |
:B...................B...B |
:B....................B..B |
:B.....................B.B |
:B......................BB |
:BBBBBBBBBBBBBBBBBBBBBBBBB |

\endsprite

\def\shboxtr{\FormOfShboxtr} 
\sprite{\FormOfShboxtr}(25,25)[0.3 em, 0.72 ex] 

:BBBBBBBBBBBBBBBBBBBBBBBBB |
:BB......................B |
:B.B.....................B |
:B..B....................B |
:B...B...................B |
:B....B..................B |
:B.....B.................B |
:B......B................B |
:B.......B...............B |
:B........B..............B |
:B.........B.............B |
:B..........B............B |
:B...........B...........B |
:B............B..........B |
:B.............B.........B |
:B..............B........B |
:B...............B.......B |
:B................B......B |
:B.................B.....B |
:B..................B....B |
:B...................B...B |
:B....................B..B |
:B.....................B.B |
:B......................BB |
:BBBBBBBBBBBBBBBBBBBBBBBBB |

\endsprite


\begin{abstract}
This is the first of two papers in which we study the modular invariance of 
pseudotraces of logarithmic intertwining operators.
We construct and study genus-one correlation functions for logarithmic
intertwining operators among generalized modules over a positive-energy and 
$C_2$-cofinite vertex operator algebra $V$.
We consider grading-restricted generalized $V$-modules which admit a right
action of some associative algebra $P$, and intertwining operators among
such modules which commute with the action of $P$ ($P$-intertwining operators).
We obtain duality properties, i.e., suitable associativity and commutativity
properties, for $P$-intertwining operators. Using pseudotraces introduced by
Miyamoto and studied by Arike, we define formal $q$-traces of products of
$P$-intertwining operators, and obtain certain identities for these formal
series. This allows us to show that the formal $q$-traces satisfy
a system of differential equations with regular singular points, and
therefore are absolutely convergent in a suitable region and can be extended
to yield multivalued analytic functions, called genus-one correlation functions.
Furthermore, we show that the space of solutions of these differential
equations is invariant under the action of the modular group.
\end{abstract}


\vspace{2em}

\section{Introduction} 

The theory of vertex operator algebra arose independently in mathematics and
physics and has been providing deep and remarkable connections between
different fields. In mathematics, one of its most spectacular applications
was the construction of the ``moonshine module'', a vertex operator
algebra (usually denoted by $V^\natural$) whose group of automorphisms
is the Monster group $\mathbb M$, the largest sporadic finite simple group.
Noticing patterns relating the dimensions of irreducible modules for
the Monster and the Fourier expansion of the modular function $J(q)$,
McKay and Thompson conjectured the existence of a ``natural'' infinite
dimensional graded module $V = \coprod_{n= -1}^\infty V_n$ for $\mathbb M$
whose graded dimension
$$
  \sum_{n=-1}^\infty \dim(V_n)q^n
$$
is given exactly by $J(q)$. Additionally, Conway and Norton conjectured
that for any element $g$ in the Monster, the series
$$
  \sum_{n\in \Z}\tr\ g|_{V_n} q^n
$$
is the Fourier expansion of a generator of the field of modular functions
for some genus zero subgroup of $SL_2(\R)$.
Frenkel, Lepowsky and Meurman constructed a module $V^\natural$ for the Monster
group in \cite{FLM}, proving the McKay-Thompson conjecture and introducing
the notion of vertex operator algebra, a variant of Borcherds' notion of
vertex algebra (\cite{B}). The full Conway-Norton conjecture
for $V^\natural$ was later proved by Borcherds.

The connection between the theory of vertex operator algebras and
the theory of modular functions has deep roots and the solution of the
Moonshine conjecture is just a part of it. 
In the important work \cite{MS1} and \cite{MS2} by Moore and Seiberg, an explicit conjecture on 
the modular invariance for intertwining operators (called "chiral vertex operators" 
in \cite{MS1} and \cite{MS2}) was stated.
In his Ph.D. thesis \cite{Z} Zhu proved a partial result on the modular invariance 
conjecture of Moore and Seiberg. Considering a class of ``rational'' vertex operator algebras satisfying
a certain cofiniteness condition, Zhu studied traces of products of $n$ vertex
operators associated to irreducible representations, and showed that these formal
traces converge and the functions thus obtained (called $n$-point genus-one
correlation functions) form a space invariant under the action of the modular
group; as direct consequence, he established the modular invariance for the
spaces of functions spanned by the graded dimension of the irreducible modules.
Zhu's results were later extended by Dong, Li and Mason in \cite{DLM2}
to include twisted representations; and in \cite{Miy1}, Miyamoto 
considered traces of products of vertex operators for modules and at most one
intertwining operator.

All these results rely heavily on the use of the commutator formula
to obtain recurrence relations for the $n$-point genus-one correlation
function in terms of the $n-1$-point functions. Since this formula
is not available for general intertwining operators, the methods do not 
generalize to product of more than one intertwining operator. In
\cite{H2}, Huang overcame this difficulty and proved the full modular invariance conjecture 
of Moore and Seiberg; he used commutativity and associativity
for intertwining operators to obtain a system of ``modular''
differential equations, and to obtain genus-one commutativity and
associativity properties.
This modular invariance result is a crucial ingredient in other important
works by Huang, including his proof of the Verlinde conjecture and the rigidity
and modularity of the vertex tensor category (see \cite{H4}, \cite{H5})
for ``rational'' vertex operator algebras.

In \cite{M1} Milas considered a class of weak modules for non-rational
vertex operator algebras, called ``logarithmic modules''. These are modules
on which the operator $L(0)$ does not necessarily act semisimply, but can
be expressed as direct sum of generalized eigenspaces for $L(0)$. Moreover,
he introduced and studied ``logarithmic intertwining operators'' between
these modules, that is, intertwining operators which involve (integral)
powers of $\log x$ in addition to powers of the formal variable $x$.
The theory of such modules and intertwining operators has
since been developed (see \cite{HLZ1}-\cite{HLZ8}) and interesting
classes of such modules have been constructed (see for instance
\cite{M2}, \cite{AM1}-\cite{AM3}).

Huang conjectured that a full modular invariance result should hold for
such classes of modules, and that it should play an important
role in the study of the properties of logarithmic modules.
Before this conjecture was explicitly formulated, a partial result
generalizing Zhu's result in the context of logarithmic modules
was obtained first by Miyamoto \cite{Miy2}, assuming only a cofiniteness
condition for the vertex operator algebra $V$ and infinite dimensionality
of all nonzero $V$-modules (an assumption which is used but not explicitly
mentioned in Miyamoto's paper, as pointed out in \cite{M2},\cite{AN}).
The main new idea is the use of a generalization of ordinary 
matrix traces, called ``pseudotraces'', to construct additional
genus-one correlation functions; pseudotraces were successively
studied by Arike in \cite{Ar}, who obtained a characterization in
terms of projective bases (or coordinate systems) for projective
modules over associative algebras. In \cite{AM4}, Milas and Adamovi\'c
considered the graded dimensions of characters of modules of certain
non-rational vertex operator (super)algebras, and proved the modularity
of the differential equations these graded dimensions satisfy.

In this paper, we obtain results that will lead us to a
full modular invariance result for logarithmic intertwining operators
generalizing the modular invariance in \cite{H2}. We consider a positive energy and 
$C_{2}$-cofinite vertex operator algebra $V$. 
We study pseudotraces
of products of intertwining operators and genus-one correlation functions
as a first step towards this result. Our first concern is on the construction of
genus-one correlation functions from products of intertwining operators;
to do so, we are naturally led to consider pseudotraces of products of
logarithmic intertwining operators. In order for the pseudotrace to be well
defined, we consider ``logarithmic'' modules which admit a right action of
some associative algebra, and logarithmic intertwining operators whose products
commute with this action. We then develop tools to study these ``formal
pseudotraces''; in particular, we formulate suitable associativity and
commutativity statements for such logarithmic intertwining operator in Theorem
\ref{assoc1} and Theorem \ref{commut1}.

Using these properties, we verify that many identities for the
formal pseudotraces in the semisimple case carry over to the logarithmic
setting; however, we see that these pseudotraces satisfy a more complicated
system of differential equations than the one in \cite{H2}
(Proposition \ref{diffeqprop}).
Nonetheless it is still possible to prove the convergence of these formal
series to multivalued analytic functions (the ``genus-one correlation
functions'') and modular invariance of the space of solutions of the system
of differential equations (Proposition \ref{diffeq-mod-inv}).

In a paper \cite{FH} in preparation jointly with Huang, we will use these results to prove that
the space of genus-one correlation functions is invariant under the
action of $SL_2(\Z)$.

\subsection{Summary of results}

In Section \ref{sec:1} we recall concepts from the theory of vertex operator algebras
and their modules. We assume that $V=\coprod_{n\in \mathbb{Z}}V_{(n)}$ 
satisfying $V_{(n)}=0$ for $n\le 0$ and $V_{(0)}=\mathbb{C}\mathbf{1}$.
We will deal with grading-restricted generalized
$V$-modules. A grading-restricted generalized
$V$-module is a direct sum of the
generalized eigenspaces for the operator $L(0)$,
$$
W = \coprod_{n\in \C} W_{[n]}
$$
equipped with a vertex operator map $Y_{W}: V\otimes W\to W((x))$
satisfying all the axioms for a $V$-module except that for $n\in \C$, the homogeneous subspaces 
$W_{[n]}$ are the
generalized eigenspaces (not necessarily eigenspaces as in the definition of $V$-module)
of $L(0) = \res_{x}xYW(\omega, x)$ with eigenvalues $n$. 
A crucial condition for obtaining differential equations for
genus-one correlation function is the $C_2$-cofiniteness condition,
introduced first in \cite{Z}:
we shall say that a grading-restricted generalized $V$-module $W$ satisfies the $C_2$-cofiniteness
condition if the space $C_2(W)$ spanned by the set
$\{v_{-2}w | v \in V, \ w\in W\}$ has finite codimension in $W$.
We assume that $V$ is $C_2$-cofinite. Then all grading-restricted generalized $V$-module 
is $C_{2}$-cofinite (see \cite{H6}).

In Section \ref{s-elliptic} we recall some notions from the theory
of elliptic functions and modular forms: in particular, we recall
the Taylor expansion of the Weierstrass elliptic function $\wp$ and
its derivatives, and the Fourier ($q$-expansion) of the Eisenstein
series. These expansions, considered as formal power series, will
appear as coefficients in identities for the formal $q$-traces 
and in the system of differential equations  for the genus-one
correlation function.

We recall the notion of pseudotrace in section \ref{sec:pt}
and define the formal $q$-traces
of products of logarithmic intertwining operators; for a fixed vertex
operator algebra $V$ and some associative algebra $P$, we consider
$V$-modules $\tilde W_i$, $i= 1,\ldots, n$ equipped with a right action
of the algebra $P$ such that $\tilde W_n$ is projective as right $P$-module; 
if $\Y_i$ are logarithmic intertwining operators of type 
$\binom{\tilde W_{i-1}}{W_i \tilde W_i}$, $i=1,\ldots, n$,
(where we take $\tilde W_0 = \tilde W_n$) such that
for all $i=1,\ldots, n$, $w_i\in W_i$, $\tilde w_i \in \tilde W_i$ and
$p\in P$,
\begin{equation*}
  \Y_i(w, x)(\tilde wp) = (\Y_i(w, x) \tilde w)p,
\end{equation*}
(we will call logarithmic intertwining operators which satisfy this property
$P$- intertwining operators) then the product
$\Y_1(w_1, x_1)\ldots \Y_n(w_n, x_n)$
is an element of
$$
  \End_P(\tilde W_n)[x_1,\dots, x_n]\{\log x_1, \dots, \log x_n\}
$$
and it is thus possible to evaluate the pseudotrace
$$
  \Tr\Y_1(w_1, x_1)\ldots \Y_n(w_n, x_n) \qL.
$$
In particular we will consider a map $\U(1): W\rightarrow W[x]$ and study
properties of formal $q$-traces obtained by taking pseudotraces of products
of \emph{geometrically modified} logarithmic intertwining operators
(\cite{H2})
$$
  \Y(\U(x)w, x),
$$
the first goal is to prove absolute convergence of such $q$-traces.

In Section \ref{s-formal-traces} we derive identities for the formal
$q$-traces; these identities have the same shape as the ones found in
\cite{H2}, and the main tools used in this section are associativity and
commutativity of intertwining operators. It is therefore necessary
to obtain a formulation of these duality properties to use in the
present context.
In Section \ref{s-duality}, under the assumptions used in \cite{HLZ7},
we state and prove suitable commutativity and associativity properties
for $P$-intertwining operators; we recall the required background from
the theory of tensor product for modules of vertex operator algebras
(\cite{HLZ1}-\cite{HLZ8}) in Section 3.2.

In Section 4 we use the identities obtained in Section \ref{s-formal-traces}
to obtain differential equation for the formal $q$-traces. We obtain a
system of differential equations for which the series
\begin{align*}
  \Tr \Y_1(\U(q_{z_1})L(0)_n^{i_1}w_1, q_{z_1})\cdot \ldots
    \cdot \Y_n(\U(q_{z_n})L(0)_n^{i_n}w_n, q_{z_n}) q^{L(0)-\frac{c}{24}},
\end{align*}
$i_j\in \N$, $j=1,\ldots,n$ are solutions 
(here $L(0)_n$ denotes the locally nilpotent part of the operator $L(0)$).
Due to the nonsemisimplicity
of the operator $L(0)$, the system is not decoupled, but nonetheless
the singular points in the variable $q$ are regular, and thus one
can prove absolute convergence of the formal $q$-traces in a suitable
domain.
We then prove modular invariance for the solutions of the system
of differential equations in Section \ref{s-diffeq-m-inv}.
We consider a space of vector valued functions in the variables
$z_1,\ldots, z_n$ and $\tau$ and we denote the components of these
vector valued functions by $\phi_{i_1,\ldots,i_n}(z_1,\ldots, z_n; \tau)$ for
$i_j\in \N$, $j=1,\ldots, n$. The solutions of the system of differential
equations are naturally elements of this space. Then, for $g\in SL_2(\Z)$,
$$
  g = \left(
    \begin{array}{c c}
      \alpha & \beta\\
      \gamma & \delta
    \end{array}
  \right)
$$
we define the action of $g$ on $\phi_{i_1,\ldots,i_n}$ by
\begin{align*}
  &(g\phi_{i_1,\ldots,i_n})(z_1,\ldots,z_n; \tau) \\
  &\qquad =\left(\frac{1}{\gamma \tau + \delta}\right)^a
    \sum_{j_1=0}^\infty \ldots \sum_{j_n=0}^\infty
    \frac{\prod_{k=1}^n(\log(\gamma\tau + \delta))^{j_k}}{j_1!\cdots j_n!}
    \phi_{i_1+j_1,\ldots,i_n+j_n}
      \left(z_1^\prime, \ldots, z_n^\prime; \tau^\prime\right)
\end{align*}
with $z^\prime = \frac{z}{\gamma \tau + \delta}$ and
$\tau^\prime = \frac{\alpha \tau + \beta}{\gamma \tau + \delta}$.
We then prove that the space of solutions of our differential
equations is invariant under this action.
\\
\\
{\bf Aknowledgements} I am grateful to Prof. Yi-Zhi Huang and Prof. 
James Lepowsky for discussions and suggestions for improvements.

\section{Logarithmic intertwining operators} \label{sec:1}

\subsection{Generalized modules and logarithmic intertwining operators}
\label{s-int-op}

After recalling some notions from logarithmic formal calculus, we recall the definitions
of vertex operator algebra and grading-restricted generalized module.  We will consider logarithmic
intertwining operators for this class of modules. For more details, see
\cite{FLM}, \cite{FHL}, \cite{LL}, \cite{H6}.

We will denote by $x$, $y$, $q$, $\log x$, $\log y$, $\log q$, $x_1$,
$x_2$, $x_3\ldots$, $\log x_1, \log x_2, \ldots$ independent
commuting formal variables.
For any set of commuting independent formal variables $X$ and
for any vector space $\mathcal W$ which does not involve any element of $X$,
we denote by $\mathcal W\{X\}$ the
space of formal series in arbitrary complex powers of the formal variables
in $X$.
In particular we will consider the space $\mathcal W \{x, \log x\}$:
an arbitrary element in this space can be written as
\begin{equation}\label{log:f}
\sum_{m,n\in {\mathbb C}}w_{n,m}x^n(\log x)^m,\;\;
w_{n,m}\in {\cal W}.
\end{equation}
The symbol $\frac{d}{dx}$ denotes the linear map (formal differentiation)
defined on $W\{x, \log x\}$ by
\begin{align}\label{ddxdef}
\lefteqn{\frac d{dx}\bigg(\sum_{m,n\in {\mathbb C}}w_{n,m}x^n(\log
x)^m\bigg)}\nn
&&\;\;\;=\sum_{m,n\in{\mathbb C}}((n+1)w_{n+1,m}+ (m+1)w_{n+1,m+1})x^{n}(\log
x)^{m}\nn
&&\bigg(=\sum_{m,n\in{\mathbb C}}nw_{n,m}x^{n-1}(\log x)^m+ \sum_{m,n\in
{\mathbb C}}mw_{n,m}x^{n-1}(\log x)^{m-1}\bigg).
\end{align}
We will make use of the notation
\begin{align*}
  \log(1 - T) &= -\sum_{n=1}^\infty \frac{T^n}{n}\\
  e^T & = \sum_{n=0}^\infty \frac{T^n}{n!},
\end{align*}
for any $T$ for which these expressions make sense.
Also, for commuting independent formal variables $x$, $y$,
we let
$$
  \log(x + y) = \log x + \log\left(1 + \frac{y}{x}\right)
    = \log x - \sum_{n=1}^\infty \frac{(-1)^n}{n}
    \left(\frac{y}{x}\right)^n.
$$

For any formal series in $W\{x, \log x\}$ the following result holds:
For $f(x)$ as in (\ref{log:f}), we have
\begin{equation}\label{log:ck1}
e^{y\frac d{dx}}f(x)=f(x+y)
\end{equation}
(``Taylor's theorem'' for logarithmic formal series) and
\begin{equation}\label{log:ck2}
e^{yx\frac d{dx}}f(x)=f(xe^y).
\end{equation}

In what follows, unless otherwise mentioned, we will fix a vertex operator algebra $V$
of cenrtal charge $c$
such that $V_{(n)} = 0$ whenever $n<0$ and $V_{(0)} = \C \one$.
We set 
$$
  C_2(V) = \text{span} \{v_{-2}u\ |\ v\in \coprod_{n>0}V_{(n)}, u\in V \}.
$$
\begin{defn} {\rm
  We say that $V$ is {\it $C_2$-cofinite} if
  $$
    \dim V/C_2(V) < \infty.
  $$
}
\end{defn}

\begin{defn}
{\rm A $\C$-graded
vector space $W=\coprod_{n\in \C}W_{[n]}$ equipped with
a linear map
\begin{eqnarray*}
Y_{W}: V\otimes W&\to &W((x))\\
v\otimes w&\mapsto & Y_{W}(v, x)w
\end{eqnarray*}
is called a {\it grading-restricted generalized $V$-module} if all the axioms for $V$-modules are satisfied except that
for $n\in \C$, the homogeneous
subspaces $W_{[n]}$ are the generalized eigenspaces of $L(0)=\res_x xY_W(\omega,x)$
with eigenvalues $n$, that is, for $n\in \C$, $w\in W_{[n]}$,
there exists $K\in \Z_{+}$, depending on $w$, such that
$(L(0)-n)^{K}w=0$. For $w\in W_{[n]}$, we denote the generalized
eigenvalue $n$ by $\wt w$.

We define {\it homomorphisms} (or {\it module maps})
and {\it isomorphisms}  between generalized $V$-modules,
{\it generalized $V$-submodules},
and {\it quotient generalized $V$-modules}
in the obvious ways.}
\end{defn}

For a grading-restricted generalized module $W$, the \emph{formal completion} of $W$ is the space 
$$\overline W = \prod_{n\in \C} W_{[n]},$$
and for $n\in \C$ we denote the projection from $\overline W$
to $W_{[n]}$ by $\pi_n$.

\begin{defn}\label{log:def}{\rm
Let $(W_1,Y_1)$, $(W_2,Y_2)$ and $(W_3,Y_3)$ be grading-restricted generalized modules
for a vertex operator algebra $V$. A {\em logarithmic
intertwining operator of type $\binom{W_3}{W_1\,W_2}$} is a linear map
\begin{equation*}
{\cal Y}(\cdot, x)\cdot: W_1\otimes W_2\to W_3[\log x]\{x\},
\end{equation*}
or equivalently,
\begin{equation*}
w_{(1)}\otimes w_{(2)}\mapsto{\cal Y}(w_{(1)},x)w_{(2)}=\sum_{n\in {\mathbb
C}}\sum_{k\in {\mathbb N}}{w_{(1)}}_{n;\,k}^{\cal Y}w_{(2)}x^{-n-1}(\log
x)^k\in W_3[\log x]\{x\}
\end{equation*}
for all $w_{(1)}\in W_1$ and $w_{(2)}\in W_2$, such that the following
conditions are satisfied: the {\em lower truncation condition}: for
any $w_{(1)}\in W_1$, $w_{(2)}\in W_2$ and $n\in {\mathbb C}$,
\begin{equation*}
{w_{(1)}}_{n+m;\,k}^{\cal Y}w_{(2)}=0\;\;\mbox{ for }\;m\in {\mathbb N}
\;\mbox{ sufficiently large,\, independently of}\;k;
\end{equation*}
the {\em Jacobi identity}:
\begin{eqnarray}\label{log:jacobi}
  \lefteqn{\dps x^{-1}_0\delta \bigg( \frac{x_1-x_2}{x_0}\bigg)
Y_3(v,x_1){\cal Y}(w_{(1)},x_2)w_{(2)}}\nn
&&\hspace{2em}- x^{-1}_0\delta \bigg( \frac{x_2-x_1}{-x_0}\bigg)
{\cal Y}(w_{(1)},x_2)Y_2(v,x_1)w_{(2)}\nn 
&&{\dps = x^{-1}_2\delta \bigg( \frac{x_1-x_0}{x_2}\bigg) {\cal
Y}(Y_1(v,x_0)w_{(1)},x_2) w_{(2)}}
\end{eqnarray}
for $v\in V$, $w_{(1)}\in W_1$ and $w_{(2)}\in W_2$ (note that the
first term on the left-hand side is meaningful because of
the lower truncation condition) and the {\em $L(-1)$-derivative property:}
for any
$w_{(1)}\in W_1$,
\begin{equation*}
{\cal Y}(L(-1)w_{(1)},x)=\frac d{dx}{\cal Y}(w_{(1)},x).
\end{equation*}
We will denote the space of all logarithmic intertwining operators
of type $\binom{W_3}{W_1W_2}$ by $\mathcal V_{W_1W_2}^{W_3}$.
}
\end{defn}
Note that if the three modules $W_1, W_2, W_3$ are ordinary modules, then
all logarithmic intertwining operators are in fact ordinary intertwining
operator, i.e., there are no logarithmic terms.

For a grading-restricted generalized $V$-module $W$, we consider the semisimple part of the
operator $L(0)$, denoted by $L(0)_s$; also denote by $L(0)_n$ the locally
nilpotent part $L(0) - L(0)_s$. Then $L(0)_n$ is a module endomorphism of
$W$ (i.e., it commutes with the action of $V$ on $W$):
\begin{lemma}
Let $\Y(x,w)$ be an intertwining operator; then
$$L(0)_n\Y(x,w) - \Y(x,w)L(0)_n = 0$$
\end{lemma}
\begin{defn} {\rm Let $W$ be a grading-restricted generalized $V$-module for a vertex operator
algebra. We define
$$ x^{\pm L(0)}: W\to W\{x\}[\log x]\subset W[\log x]\{x\}$$
by
\begin{equation*}
x^{\pm L(0)} = x^{\pm L(0)_s}e^{\pm \log x (L(0)-L(0)_s)}
\end{equation*}
}
\end{defn}

\begin{lemma}\label{x-L-0-derivative}
  Using the same notation as above, we have
  \begin{equation*}
    \frac{d}{dx} x^{L(0)} = L(0)x^{L(0)-1}.
  \end{equation*}
\end{lemma}

\begin{prop}
  Let ${\cal Y}$ be a logarithmic intertwining operator of type
  $\binom{W_3}{W_1\,W_2}$ and let $w\in W_1$. Then
  \begin{description}
    \item{(a)}
    \begin{equation*}
      e^{yL(-1)}{\cal Y}(w,x)e^{-yL(-1)}={\cal Y}(e^{yL(-1)}w,x)={\cal Y}(w,x+y)
    \end{equation*}
    \item{(b)}
    \begin{equation*}
      y^{L(0)}{\cal Y}(w,x)y^{-L(0)}={\cal Y}(y^{L(0)}w,xy)
    \end{equation*}
    \item{(c)}
    \begin{equation*}
      e^{yL(1)}{\cal Y}(w,x)e^{-yL(1)}={\cal
      Y}(e^{y(1-yx)L(1)}(1-yx)^{-2L(0)}w,x(1-yx)^{-1}).
    \end{equation*}
  \end{description}
\end{prop}

\subsection{Tensor product of modules and associativity of logarithmic intertwining operators}
\label{s-tensor}

In this section, we recall the notion of $P(z)$-tensor product for
modules of vertex operator algebras introduced in \cite{HL1}-\cite{HL3},
\cite{H3} (and then generalized to the logarithmic
case  in \cite{HLZ1}-\cite{HLZ8} and \cite{H6}), and some related results that we will need
later; in particular, we will state associativity for logarithmic intertwining operators.

Let $W_1, W_2, W_3$ be grading-restricted generalized $V$-modules. 
Related to the concept of logarithmic intertwining operator is that of \emph{intertwining map}:
a $P(z)$-intertwining map of type $\binom{W_3}{W_1W_2}$ is a linear function
\begin{align*}
  I_z: W_1\otimes W_2 \rightarrow \overline W_3
\end{align*}
satisfying the following conditions: the
\emph{lower truncation condition}: for $w_{(1)} \in W_1$ and $w_{(2)}\in W_2$
and any $n\in \C$,
$$\pi_{n-m}I(w_{(1)}\otimes w_{(2)}) = 0 \text{ for }m\in \N 
  \text{ sufficiently large;}
$$
and the \emph{Jacobi identity}: for $v\in V$, $w_{(1)}\in W_1$ and $w_{(2)}\in W_2$,
\begin{align*}
  &\dlt{x_0}{x_1}{-z}Y_3(v, x_1)I(w_{(1)}\otimes w_{(2)})\\
    &\qquad = \dlt{z}{x_1}{-x_0}I(Y_1(v, x_0)w_{(1)}\otimes w_{(2)})\\
      &\qquad \qquad +\dlt{x_0}{-z}{+x_1}I(w_{(1)}\otimes Y_2(v,x_1)w_{(2)}).
\end{align*}
Let $\Y\in \binom{W_3}{W_1W_2}$, and let $p\in \Z$; then the map
\begin{align*}
  I_{\Y,p} : W_1\otimes W_2 &\rightarrow W_3\\
    w_{(1)}\otimes w_{(2)} &\mapsto \Y(w_{(1)}, x)w_{(2)} 
  \Big|_{x^n=e^{n(\log z + 2\pi i p)}, (\log(x))^m = (\log z + 2\pi i p)^m}
\end{align*}
is a well defined $P(z)$-intertwining map of type $\binom{W_3}{W_1W_2}$; for
any $p\in \Z$, the correspondence $\Y \mapsto I_{\Y,p}$ is a bijection with
inverse denoted by $I\mapsto \Y_{I,p}$ (\cite{HLZ3}).

\begin{defn}
  [\cite{HLZ1}-\cite{HLZ7}] {\rm Given two grading-restricted generalized $V$-modules $W_1$,
  $W_2$, their $P(z)$-tensor product is a third such module, denoted by
  $W_1\boxtimes_{P(z)} W_2$, equipped with a $P(z)$-intertwining map
  \begin{equation*}
    \boxtimes_{P(z)} : W_1\otimes W_2 \rightarrow
      \overline {W_1\boxtimes_{P(z)} W_2}
  \end{equation*}
  such that for any grading-restricted generalized $V$-module $W_3$ and intertwining map
  $I: W_1\otimes W_2 \rightarrow \overline W_3$, there exists a unique
  $V$-module morphism $\eta : W_1\boxtimes_{P(z)}W_2 \rightarrow W_3$
  such that
  \begin{equation*}
    I = \overline\eta\circ \boxtimes_{P(z)},
  \end{equation*}
  where $\overline \eta$ is the unique map $\overline \eta :
    \overline {W_1\boxtimes_{P(z)} W_2} \rightarrow \overline W_3$
  extending $\eta$.
  }
\end{defn}
From the definition, one can see that given two $V$-modules, if their
$P(z)$-tensor product exists then it is unique; moreover we have the
following:
\begin{prop}\label{tensor-morphism}
  Let $W_1$, $W_2$, $W_3$ and $W_4$ be generalized $V$-modules and
  $\varphi:W_1\rightarrow W_3$ and $\psi:W_2\rightarrow W_4$
  be $V$-module homomorphisms.
  Suppose that the $P(z)$-tensor products $W_1\boxtimes_{P(z)}W_2$ and
  $W_3\boxtimes_{P(z)}W_4$ exist and denote their intertwining maps by
  $I_1$, $I_2$ respectively.
  Then there exists a unique $V$ homomorphism
  \begin{equation*}
    \varphi\boxtimes_{P(z)}\psi : W_1\boxtimes_{P(z)} W_2 \rightarrow W_3\boxtimes_{P(z)} W_4
  \end{equation*}
  such that for all $w_{(1)}\in W_1$ and $w_{(2)}\in W_2$,
  \begin{equation*}
    I_2(\varphi(w_{(1)})\otimes\psi(w_{(2)})) =
      \overline {\varphi\boxtimes_{P(z)}\psi}\circ I_1(w_{(1)}\otimes w_{(2)}).
  \end{equation*}
\end{prop}
We recall some results from \cite{H6} and \cite{HLZ5} concerning the existence of $P(z)$-tensor product and
associativity of logarithmic
intertwining operators.

\begin{thm} [\cite{HLZ5}] \label{HLZassociativity}
  Let $V$ be a vertex operator algebra whose category of grading-restricted generalized modules is
  closed under $P(z)$-tensor product, and such that every grading-restricted generalized module satisfies
  the $C_1$-cofiniteness condition and the quasi-finite dimensionality condition;
  let $W_1$, $W_2$, $W_3$, $W_4$, $M$ be grading-restricted generalized $V$-modules.
  \begin{enumerate}
    \item Consider intertwining operators $\Y_1\in \binom{W_4}{W_1M}$,
      $\Y_2\in \binom{M}{W_2W_3}$. Then there exists a unique intertwining operator
      $\Y^1 \in \binom{W_4}{W_1\boxtimes_{P(z_0)} W_2\ W_3}$ such that
      \begin{align*}
        &\langle w_{(4)}^\prime, \Y_1(w_{(1)}, x_1)
          \Y_2(w_{(2)},x_2)w_{(3)}\rangle\Big|_{x_i^n=e^{n \log z_i},\; (\log x_i)^m = (\log z_i)^m,\; i=1,2}\\
        &= \langle w_{(4)}^\prime,
          \Y^1(\Y_{\boxtimes_{P(z_0),0}}(w_{(1)},x_0)w_{(2)},x_2)
            w_{(3)}\rangle\Big|_{x_i^n=e^{n \log z_i},\; (\log x_i)^m = (\log z_i)^m,\; i=0,2}
        \end{align*}
    whenever $z_0 = z_1 - z_2$ and $|z_1| > |z_2| > |z_0| > 0$, for all
      $w_{(1)}\in W_1, w_{(2)}\in W_2,w_{(3)}\in W_3$ and $w_{(4)}^\prime \in W_4^\prime$.
    \item Let $\Y^1 \in \binom{W_4}{M W_3}$, $\Y^2 \in \binom{M}{W_1 W_2}$. Then
      there exists a unique intertwining operator $\Y_1$ of type
      $\binom{W_4}{W_1 W_2\boxtimes_{P(z_2)}W_3}$ such that
      \begin{align*}
        &\langle w_{(4)}^\prime, \Y^1(\Y^2(w_{(1)}, x_0)w_{(2)}, x_2)
          w_{(3)}\rangle\Big|_{x_i^n=e^{n \log z_i},\; (\log x_i)^m = (\log z_i)^m,\; i=0,2}\\
        &\qquad = \langle w_{(4)}^\prime, \Y_1(w_{(1)},x_1)
          \Y_{\boxtimes_{P(z_2),0}}(w_{(2)},x_2)w_{(3)}\rangle\Big|_{x_i^n=e^{n \log z_i},\; (\log x_i)^m = (\log z_i)^m,\; i=1,2}
      \end{align*}
      whenever $|z_1| > |z_2| > |z_1 - z_2| > 0$ for all
      $w_{(1)}\in W_1, w_{(2)}\in W_2,w_{(3)}\in W_3$ and $w_{(4)}^\prime \in W_4^\prime$.
  \end{enumerate}
\end{thm}

\begin{thm}[\cite{H6}]
 Assume that $V$ is $C_2$-cofinite. Then the category of grading-restricted $V$-modules is closed 
under $P(z)$-tensor product and every grading-restricted $V$-module is $C_2$-cofinite.
\end{thm}

\subsection{Elliptic functions and Eisenstein series}
\label{s-elliptic}

In this section we recall some basic properties of the Weierstrass $\wp$ function
and Eisenstein series; in particular, we will be using the
Taylor and Fourier $q$-expansions of such functions. For additional
background, see \cite{L,Z}.
For $z\in \C$, we will use the notation $q_z = e^{2\pi i z}$.
We first introduce a formal power series related to the expansion
of the Weierstrass $\wp$ function:
for $m\geq 0$,
\begin{align*}
  P_{m+1}(x;q) = (2\pi i)^{m+1}\sum_{l>0}
  \left(\frac{l^m}{m!}\frac{x^l}{1-q^l} - 
  \frac{(-1)^ml^m}{m!}\frac{q^lx^{-l}}{1-q^l}\right)
\end{align*}
where $(1-q^l)^{-1}$ is the power series $\sum_{k\geq 0} q^{lk}$ in the formal variable
$q$. For $\tau, z\in \C$ satisfying $|q_\tau| < |q_z| < 1$, the series
$P_{m+1}(q_z;q_\tau)$ is absolutely convergent, and for $|q_z|<1$ the
$q$-coefficients of $P_{m+1}(q_z;q)$ are absolutely convergent.
Let
\begin{align*}
  \wp_1(z;\tau) &= \frac{1}{z} + \sum_{(k,l)\neq(0,0)}
      \left(\frac{1}{z-(k\tau+l)}
      + \frac{1}{k\tau + l}
      + \frac{z}{(k\tau + l)^2}\right)\\
  \wp_2(z;\tau) &= \frac{1}{z^2} + \sum_{(k,l)\neq(0,0)}
      \left(\frac{1}{(z-(k\tau+l))^2}
      - \frac{1}{(k\tau + l)^2}\right);
\end{align*}
and for $m \geq 2$, let
\begin{align*}
  \wp_{m+1}(z;\tau) = -\frac{1}{m}\frac{\partial}{\partial z} \wp_m(z;\tau).
\end{align*}
These functions have Laurent expansion
\begin{align*}
  \wp_m(z;\tau) = \frac{1}{z^m} + (-1)^m\sum_{k\geq 1}
    \binom{2k+1}{m-1} G_{2k+2}(\tau)z^{2k+2-m}
\end{align*}
in the region $0<|z|< \min(1, |\tau|)$, where $G_{2k+2}(\tau)$ are the Eisenstein
series defined by
\begin{align*}
G_{2k+2}(\tau) = \sum_{(m,l)\in\Z^2\setminus (0,0)} \frac{1}{(m\tau + l)^{2k+2}}
\end{align*}
for $k\geq 1$. Moreover, let
\begin{align*}
  G_2(\tau) = \frac{\pi^2}{3} + \sum_{m\in \Z\setminus \{0\}}
    \sum_{l\in \Z} \frac{1}{(m\tau + l)^2}.
\end{align*}
It is known that the Eisenstein series have $q$-expansion
\begin{align*}
  G_{2k+2}(\tau) = 2\zeta(2k+2) + \frac{2(2\pi i)^{2k+2}}{(2k+1)!}
  \sum_{l=1}^\infty \frac{l^{2k+1}q_\tau^l}{1 - q_\tau^l}
\end{align*}
We use the following notation to denote these as formal power series in the
variable $q$:
\begin{align*}
  \tilde G_{2k+2}(q) = 2\zeta(2k+2) + \frac{2(2\pi i)^{2k+2}}{(2k+1)!}
  \sum_{l=1}^\infty \frac{l^{2k+1}q^l}{1 - q^l},\qquad k \in \N
\end{align*}
and similarly for the expansion of the elliptic functions
\begin{align*}
  \tilde \wp_m(x; q) = \frac{1}{x^m} + (-1)^m \sum_{k=1}^\infty
    \binom{2k+1}{m-1} \tilde G_{2k+2}(q) x^{2k+2 - m},\qquad m=1,2,\ldots
\end{align*}
For $z\in \C$ such that $0 < |z| < 1$,
\begin{align*}
  \tilde \wp_m(z;q) = (-1)^m
    \left(P_m(q_z;q) - \der{m-1}{z}(\tilde G_2(q)z - \pi i)\right)
\end{align*}
and $\tilde \wp_m(z; q_\tau) = \wp_m(z; \tau)$. The following is well known:
\begin{prop}
  For any element $g= \left(
      \begin{array}{r r}
        \alpha & \beta\\
        \gamma & \delta
      \end{array}
    \right)\in SL_2(\Z)$, and $m=1,2,\ldots$
  \begin{align*}
    \wp_m|_g(z; \tau) := (\gamma \tau + \delta)^{-m}
      \wp_m\left(\frac{z}{\gamma \tau + \delta};
      \frac{\alpha\tau + \beta}{\gamma\tau+\delta}\right) = \wp_m(z;\tau)
  \end{align*}
  if $m > 1$,
  \begin{align*}
    \wp_m(z+\tau; \tau) = \wp_m(z+1; \tau) = \wp_m(z; \tau)
  \end{align*}
  and
  \begin{gather*}
    \wp_1(z+1;\tau) = \wp_1(z;\tau) + G_2(q)\\
    \wp_1(z+\tau; \tau) = \wp_1(z;\tau) + G_2(\tau)\tau -2\pi i
  \end{gather*}
\end{prop}

\begin{prop}
  Let $f(q)$ be a modular form of weight $k$. Then the function $\vartheta_k(f)$
  defined by
  \begin{equation*}
    (2\pi i)^2 q \frac{d}{dq}f(q) + k G_2(q)f(q)
  \end{equation*}
  is a modular form of weight $k+2$.
\end{prop}

\subsection{Pseudotraces}\label{sec:pt}

Let $P$ be a finite-dimensional associative algebra over $\C$.
We will say that a linear function $\phi: P\rightarrow \C$ is \emph{symmetric}
if  $\phi(pq) = \phi(qp)$ for all $p,q\in P$, and denote by $SLF(P)$ the
vector space of all such functions; $SLF(P)\simeq \left(P/[P,P]\right)^*$.

Let $M$ be a finitely generated projective right $P$-module. 
It is well known that this condition is equivalent to the existence of a
\emph{projective basis} for $M$, that is
a pair of sets $\{m_i\}_{i=1}^n\subseteq M$,
$\{\alpha_i\}_{i=1}^n\subseteq \hom_P(M, P)$
such that for all $m\in M$, $$m = \sum_{i=1}^n m_i\alpha_i(m).$$ 

\begin{defn}{\rm Let $\phi \in SLF(P)$.
    The \emph{pseudtrace map} $\phi_M$ on $\End_P(M)$ associated to $\phi$
    is the function $\phi_M$ defined by
\begin{align*}
  \phi_M(\alpha) = \phi\left(\sum_{i=1}^n \alpha_i(\alpha(m_i))\right)
\end{align*}
for all $\alpha \in \End_P(M)$.
}
\end{defn}
It is easy to show that the definition of pseudotrace does not depend
on the choice of projective basis for $M$.
Pseudotrace maps are an extension of regular trace functions, and share similar
properties; it is easy to see that the pseudotrace of a product is invariant under
cyclic permutation of the factors \cite{AN}:

\begin{prop}
  Let $M_1$ and $M_2$ be two right projective $P$-modules, and consider
  homomorphisms
  $\alpha \in \hom_P(M_1,M_2)$ and $\beta \in \hom_P(M_2,M_1)$; then
  $$
    \phi_{M_1} (\beta \circ \alpha) =
    \phi_{M_2} (\alpha \circ \beta).
  $$
\end{prop}

Importantly, for a given associative algebra $P$, the pseudotraces of
representations of $P$ span the space $SLF(P)$ of symmetric linear functions
on $P$ (see for example \cite{Ar, Miy2}).

\subsection{Formal $q$-traces of logarithmic intertwining operators}
\label{s-formal-traces}

In this section we consider grading-restricted $V$-modules
which admit a right action (by module endomorphisms) of a finite-dimensional associative
algebra $P$. We then consider products of intertwining operators which
commute with this action; in particular, we are able to define the formal
$q$-trace of such products by using pseudotraces on the $P$-endomorphism
ring of the $L(0)$ generalized eigenspaces in the grading-restricted generalized $V$-modules.
We use properties of pseudotraces and intertwining operators to derive
identities for the formal $q$-traces.

Let $P$ be a finite-dimensional associative algebra equipped with a symmetric linear function
$\phi$.
We say that a grading-restricted generalized $V$-module $W$ is a \emph{$V$-$P$-bimodule} if $W$
is a right $P$ module and $P$ acts on $W$ by $V$-module endomorphisms, that
is, for any $v\in V$, $w\in W$, and $p \in P$,
$$
  Y(v,x)(wp) = (Y(v,x)w)p.
$$
\begin{prop} \label{Wn-P-submodule}
  Let $W_{[n]}$ be the $L(0)$-generalized eigenspace of $W$ for the eigenvalue
  $n$. Then $W_{[n]}$ is a $P$-submodule of $W$; if $W$ is a projective right
  $P$-module, so is $W_{[n]}$.
\end{prop}
\proof
  This is clear since the action of $P$ commutes with $L(0)$:
  let $p\in P$ and $w\in W_{[n]}$; then there exists $k\in \N$ such that
  $(L(0) - n)^kw = 0$. Then
  $$
  (L(0)-n)^k(wp) = ((L(0)-n)^kw)p = 0
  $$
  which proves the first part of the claim; the second part follows since
  $W_{[n]}$ is a direct summand of $W$.
\epfv

\noindent As a consequence, the action of $P$ commutes with $L(0)_s$ and $L(0)_n$.
Suppose now $W$ is a grading-restricted generalized $V$-module,
projective as right $P$-module. Then for any generalized
eigenspace $W_{[n]}$, we have the pseudotrace 
$$
  \phi_{W_{[n]}} : \End_P(W_{[n]}) \rightarrow \C;
$$
and for a given element
$a(x_1,\ldots,x_k) \in \End_P(W)\{x_1,\ldots,x_k, \log x_1,\ldots, \log x_k\}$,
we define
$$
  \tr^\phi_W a(x_1,\ldots,x_n)q^{L(0)} =
    \sum_{n\in \C} \phi_{W_{[n]}} \left.\left(\pi_n a(x_1,\ldots,x_k)
    \sum_{i=0}^\infty\frac{(L(0)_n)^i}{i!}(\log q)^i\right)\right|_{{W_{[n]}}}q^n
$$
where $\pi_n : W \rightarrow W_{[n]}$ is the projection on the generalized eigenspace
$W_{[n]}$. Note that since $L(0)_n$ is locally nilpotent, the summation over $i$ is
finite for any value of $n\in \C$. If $W$ has finite length $l$, then the powers
of $\log q$ are globally bounded by $l$.
\begin{rema} \label{cfremark}
  {\rm
    Suppose $T\in \End_P \left(W_{[n]}\right)$, and let $\{w_i\}_{i=1}^s$,
    $\{\alpha_i\}_{i=1}^s$ be a projective basis for $W_{[n]}$.
    Let $w_i^\prime \in W_{[n]}^*$ be the linear function defined by
    $\langle w_i^\prime, w\rangle = \phi(\alpha_i(w))$ for all $w\in W_{[n]}$.
    Then one can express the pseudotrace of $T$ as
    \begin{equation*}
      \phi_{W_{[n]}}(T) = \sum_{i=1}^s \langle w_i^\prime, T w_i\rangle.
    \end{equation*}
    In particular, for any $n\in \C$ let $\{w_{n,i}\}_{i=1}^{s_n}$,
    $\{\alpha_{n,i}\}_{i=1}^{s_n}$ be a projective basis of $W_{[n]}$
    and let $a(x_1,\ldots, x_k)$ as above, one can express the $q$-trace
    of $a(x_1, \ldots, x_n)$ as
    \begin{equation*}
      \tr_W^\phi a(x_1,\ldots, x_n) q^{L(0)}=
        \sum_{n\in \C} \sum_{i=1}^{s_n}\left( \sum_{j=1}^{\infty} 
          \left\langle w_{n,i}^\prime, a(x_1,\ldots, a_n)
            \frac{(L(0)_n)^j}{j!} w_{n,i}\right\rangle (\log q)^j \right) q^n
    \end{equation*}
    where $w_{n,i}^\prime$ is defined as above and extended to an element
    of $W^\prime$ by letting it map to $0$ the generalized subspaces for eigenvalues
    different from $n$.
  }
\end{rema}

Now let $W_i$, $\tilde W_i$, $i=1,\ldots,n$, be grading-restricted generalized modules
for $V$, and suppose $\tilde W_n$ is a $V$-$P$-bimodules, projective as right
$P$-module.
Moreover, consider logarithmic intertwining operators 
$\Y_i \in \binom{\tilde W_{i-1}}{W_i\tilde W_i}$, $i=1,\ldots,n$, where we use the
convention $\tilde W_0 = \tilde W_n$. If the action of $P$ commutes with the
product of these intertwining operators, i.e., for all $w_i\in W_i$,
$\tilde w_n \in \tilde W_n$, and $p$ in $P$,
$$
  \Y_1(w_1, x_1)\cdots\Y_n(w_n, x_n)(\tilde w_np)
    = (\Y_1(w_1, x_1)\cdots\Y_n(w_n, x_n)\tilde w_n)p,
$$
then we have the \emph{formal $q$-trace}
$$
  \Tr \Y_1(w_1, x_1)\cdots\Y_n(w_n, x_n)q^{L(0)}.
$$

\begin{rema}{\rm
  Using the same notation as above, suppose  $\tilde W_i$ for $i=1,\ldots,n$
  are $V$-$P$-bimodules, with $\tilde{W}_n$ projective as right $P$-module, and the action
  of $P$ commutes with all the intertwining operators, i.e., for all $i=1,\ldots,n$
  for all $w_i\in W_i$, $\tilde w_i \in \tilde W_i$, and $p$ in $P$,
  $$
    \Y(w_i, x) (\tilde w_ip) = (\Y(w_i, x)\tilde w_i)p.
  $$
  Then the product of the intertwining operators commute with the action of $P$ and
  $$
    \Tr \Y_1(w_1, x_1)\cdots\Y_n(w_n, x_n)q^{L(0)}
  $$
  is well defined.}
\end{rema}

\begin{defn}{\rm
  Let $W$ be a grading-restricted generalized $V$-module and $W_1, W_2$ be $V$-$P$-bimodules; we will say that
  a logarithmic intertwining operator $\Y$ of type $\binom{W_2}{W W_1}$ is a
  \emph{$P$-intertwining operator} if
  $$
  \Y(w, x)(w_{1}p) = (\Y(w, x)w_1)p
  $$
  for all $p\in P$, $w\in W$, $W_1\in W_1$.
  }
\end{defn}

\begin{rema} {\rm
  Note that for any set of intertwining operators $\Y_1,\ldots,\Y_n$ of
  the above types, one can always consider $\tilde W_i$ as a projective
  right $P$-module with $P=\C$ and $\phi = 1_{\C}$. In this case, $\tr^\phi$
  corresponds to the ordinary matrix trace of the action of the intertwining
  operators on the vector space $\tilde W_n$.
}
\end{rema}

\begin{rema}
  {\rm
  Using the same notation as in Remark \ref{cfremark}, we can express the formal
  $q$-trace as
  \begin{align*}
    &\Tr \Y_1(w_1, x_1)\cdots\Y_n(w_n, x_n)q^{L(0)} \nn
    &\qquad = \sum_{n\in \C} \sum_{i=1}^{s_n}\left( \sum_{j=1}^{\infty} 
          \left\langle w_{n,i}^\prime, 
            \Y_1(w_1, x_1)\cdots\Y_n(w_n, x_n)
            \frac{(L(0)_n)^j}{j!} w_{n,i}\right\rangle (\log q)^j \right) q^n.
  \end{align*}
  In particular, as a formal series in the variables $q$, $\log q$,
  its coefficients are finite sums of genus-zero correlation functions; one
  will be able to use properties of these correlation functions to obtain
  properties for formal $q$ traces.
  }
\end{rema}

Following \cite{H2}, we use the concept of \emph{geometrically modified
intertwining operator}. Let $A_j$, $j\in \Z_+$ be the numbers defined by
the formal relation
$$
  \frac{1}{2\pi i} \log(1+2\pi i y) = \left( 
  \exp \left( \sum_{j\in \Z_+} A_j y^{j+1}\frac{\partial}{\partial y} \right)
  \right)y
$$
and let $L_+(A) = \sum_{j\in Z_+}A_j L(j)$;
then the operator $\U(1)$ is defined by
$$
  \U(1) = (2\pi i)^{L(0)}e^{-L_+(A)}.
$$
Also let $\U(x) = x^{L(0)}\U(1)$, for any formal expression $x$ for which the
expression makes sense: it follows that
$$
y^{L(0)}\U(x) = \U(yx).
$$
\begin{defn}[\cite{H2}] {\rm Let $W_1$, $W_2$, $W_3$ be grading-restriced generalized
$V$-modules, and $\Y$ a
logarithmic intertwining operator of type $\binom{W_3}{W_1 W_2}$.
The operator $\Y(\U(x)w, x)$ is called a \emph{geometrically modified
logarithmic intertwining operator}. }
\end{defn}
Here we recall some of the properties of the operator $\U(1)$ and
 the geometrically modified logarithmic intertwining operators; see \cite{H2} for
these results (the proofs in \cite{H2} are for ordinary intertwining operators
but the proofs for logarithmic intertwining operators are exactly the same).
\begin{lemma} 
  Let $\Y$ be a logarithmic intertwining operator of type $\binom{W_3}{W_1,W_2}$ for grading
  restricted generalized $V$-modules $W_1$, $W_2$, $W_3$.
  Then for $u\in V$ and $w\in W_1$,
  \begin{align} \label{geom-commutator}
    &[Y(\U(x_1)u,x_1),\Y(\U(x_2)w,x_2)] \nn
    &\qquad= 2\pi i \res_y \delta \left( \frac{x_1}{e^{2\pi iy}x_2} \right)
    \Y(\U(x_2)Y(u,y)w,x_2)
  \end{align}
\end{lemma}
\begin{lemma}
  Let $W_1$, $W_2$, $W_3$ be grading-restricted generalized $V$-modules,
  and $\Y$ and logarithmic intertwining operator of type
  $\binom{W_3}{W_1 W_2}$.
  Then for any $w_1\in W_1$,
  \begin{equation} \label{L-1-derivative}
    \Y(\U(x)L(-1)w_1,x) = 2\pi ix \frac{d}{dx}\Y(\U(x)w_1, x).
  \end{equation}
\end{lemma}
\begin{lemma}
  For any grading-restricted generalized $V$-module $W$ and $u\in V$, we have
  the $\U(x)$ conjugation property
  \begin{equation} \label{U-conj}
    \U(x)Y(u,y) = Y(\U(xe^{2\pi i y})u, x(e^{2\pi i y} -1))\U(x).
  \end{equation}
\end{lemma}
We will consider formal $q$-traces of products of geometrically modified logarithmic
intertwining operators
$$
  \Tr \Y_1(\U(x_1)w_1, x_1)\cdots\Y_n(\U(x_n)w_n, x_n)\qL
$$
for intertwining operators whose product commutes with the action of $P$.
Many of the properties of regular traces which hold in the completely reducible case
carry over to the logarithmic setting.

In the following, for any $v\in V$, we denote by $o(v)$ the constant term of
the operator $Y(x^{L(0)}v, x)$ acting on a generalized $V$-module; that is,
$$
  o(v) = \res_x x^{-1}Y(x^{L(0)}v,x) = v_{\swt v - 1}.
$$

\begin{lemma}
  Consider grading-restricted generalized $V$-modules $W_i$, $\tilde W_i$
  for $i=1,\ldots,n$, with $\tilde W_0 = \tilde W_n$, and logarithmic intertwining
  operators $\Y_i$ of type $\binom{\tilde W_{i-1}}{W_i \tilde W_i}$ for
  $i=1,\ldots,n$. Moreover, suppose $\tilde W_0$ is a $V$-$P$-bimodule
  projective as a right $P$-module for some algebra
  $P$ equipped with a symmetric linear function $\phi$, and that the product
  of the intertwining operators $\Y_1,\ldots, \Y_n$ commutes with the action
  of $P$.
  Then for any $v\in V$, $w_i\in W_i$, we have
  \begin{align}\label{identity1}
    &\Tr Y(\U(x)u,x)\Y_1(\U(x_1)w_1,x_1)\cdots
    \Y_n(\U(x_n)w_n,x_n)q^{L(0)}\nn
    &\qquad = \sum_{i=1}^n\sum_{m\geq 0} P_{m+1}\left(\frac{x_i}{x}; q\right)
    \Tr \Y_1(\U(x_1)w_1,x_1)\cdot\nn
    &\qquad\qquad\qquad\cdots \Y_{i-1}(\U(x_{i-1})w_{i-1},x_{i-1})
      \Y_i(\U(x_i)u_m w_i,x_i)\cdot\nn
    &\qquad\qquad\qquad\cdot \Y_{i+1}(\U(x_{i+1})w_{i+1},x_{i+1})\cdots
      \Y_{n}(\U(x_{n})w_{n},x_{n})q^{L(0)}\nn
    &\qquad + \Tr o(\U(1)u)\gY 1\cdots \gY n \qL
  \end{align}
  and
  \begin{align}\label{identity2}
    &\sum_{i=1}^n \Tr \gY 1\cdots
      \gY {i-1}\cdot\nn
      &\qquad\cdot \Y_i(\U(x_i)u_0w_i,x_i)
      \gY{i+1} \cdots\gY n q^{L(0)}\nn
      &\qquad = 0
    \end{align}
\end{lemma}

\proof
  By the commutator formula,
  \begin{align*}
    &Y(\U(x)u,x)\gY 1 \cdots\gY n\\
    &\quad = \sum_{i=1}^n \Y_1(\U(x_1)w_1,x_1)
      \cdots\gY{i-1}\cdot\\
    &\qquad\qquad \cdot [Y(\U(x)u,x), \gY i]
      \gY {i+1}\cdot\\
    &\qquad\qquad \cdots \gY n\\
    &\quad\qquad + \gY 1\cdots\gY n Y(\U(x)u, x)\bn
    &\quad = \sum_{i=1}^n 2\pi i \res_y \delta\left(\frac{x}{e^{2\pi iy}x_i}\right)
      \Y_1(\U(x_1)w_1,x_1)
      \cdots\gY{i-1}\cdot\\
    &\qquad\qquad \cdot \Y_i(\U(x_1)Y(u,y)w_i,x_i)
      \gY {i+1}\cdot\\
    &\qquad\qquad \cdots \gY n\\
    &\quad\qquad + \gY 1\cdots\gY n Y(\U(x)u, x).
  \end{align*}
  Since $P$ acts on $\tilde W_n$ as $V$-module endomorphisms,
  and the product of the intertwining operators $\Y_1\ldots\Y_n$ commutes
  with $P$, the pseudotrace of each term in this expression is well defined.
  Therefore, by linearity of pseudotraces, $\qL$ conjugation property and cyclic
  property of pseudotraces,
  \begin{align*}
    &\Tr Y(\U(x)u,x)\gY 1 \cdots\gY n \qL\\
    &\quad = \sum_{i=1}^n \Tr 2\pi i \res_y \delta\left(\frac{x}{e^{2\pi iy}x_i}\right)
      \Y_1(\U(x_1)w_1,x_1)
      \cdots\gY{i-1}\cdot\\
    &\qquad\qquad \cdot \Y_i(\U(x_1)Y(u,y)w_i,x_i)
      \gY {i+1}\cdot\\
    &\qquad\qquad \cdots \gY n\qL\\
    &\quad\qquad + \Tr \gY 1\cdots\gY n Y(\U(x)u, x)\qL \bn 
    &\quad = \sum_{i=1}^n \Tr 2\pi i \res_y \delta\left(\frac{x}{e^{2\pi iy}x_i}\right)
      \Y_1(\U(x_1)w_1,x_1)
      \cdots\gY{i-1}\cdot\\
    \displaybreak[1]
    &\qquad\qquad \cdot \Y_i(\U(x_1)Y(u,y)w_i,x_i)
      \gY {i+1}\cdot\\
    &\qquad\qquad \cdots \gY n\qL\\
    \displaybreak[1]
    &\quad\qquad + \Tr \gY 1\cdots\gY n \qL 
      Y\left(\U\left(\frac x q\right)u,\frac x q\right)\bn
    &\quad = \sum_{i=1}^n \Tr 2\pi i \res_y \delta\left(\frac{x}{e^{2\pi iy}x_i}\right)
      \Y_1(\U(x_1)w_1,x_1)
      \cdots\gY{i-1}\cdot\\
    &\qquad\qquad \cdot \Y_i(\U(x_1)Y(u,y)w_i,x_i)
      \gY {i+1}\cdot\\
    &\qquad\qquad \cdots \gY n\qL\\
    &\quad\qquad + \Tr Y\left(\U\left(\frac x q\right)u,\frac x q\right)
      \gY 1\cdots\gY n \qL\bn
    &\quad = \sum_{i=1}^n \Tr 2\pi i \res_y \delta\left(\frac{x}{e^{2\pi iy}x_i}\right)
      \Y_1(\U(x_1)w_1,x_1)
      \cdots\gY{i-1}\cdot\\
    &\qquad\qquad \cdot \Y_i(\U(x_1)Y(u,y)w_i,x_i)
      \gY {i+1}\cdot\\
    &\qquad\qquad \cdots \gY n\qL\\
    &\quad\qquad + \left(q^{-x\frac{\partial}{\partial x}}\right)\Tr Y(\U(x)u,x)
      \gY 1\cdots\gY n \qL 
  \end{align*}
  and thus
  \begin{align*}
    &\left(1-q^{-x\frac{\partial}{\partial x}}\right)\Tr Y(\U(x)u,x)
      \gY 1\cdots\gY n \qL\bn
    &\quad = \sum_{i=1}^n \Tr 2\pi i \res_y \delta\left(\frac{x}{e^{2\pi iy}x_i}\right)
      \Y_1(\U(x_1)w_1,x_1)
      \cdots\gY{i-1}\cdot\\
    &\qquad\qquad \cdot \Y_i(\U(x_1)Y(u,y)w_i,x_i)
      \gY {i+1}\cdot\bn
    &\qquad\qquad \cdots \gY n\qL\bn
    &\quad = \sum_{i=1}^n \Tr 2\pi i \res_y e^{2\pi i yx_i\frac{\partial}{\partial x_i}}
      \delta\left(\frac{x}{x_i}\right)
      \Y_1(\U(x_1)w_1,x_1)
      \cdots\gY{i-1}\cdot\\
    &\qquad\qquad \cdot \Y_i(\U(x_1)Y(u,y)w_i,x_i)
      \gY {i+1}\cdot\\
    &\qquad\qquad \cdots \gY n\qL\bn
    &\quad = \sum_{i=1}^n\sum_{m=0}^\infty
      \Tr \frac{(2\pi i)^{m+1}}{m!} \left(x_i\frac{\partial}{\partial x_i}\right)^m
      \delta\left(\frac{x}{x_i}\right)\cdot\bn
    &\qquad\qquad \cdot\Y_1(\U(x_1)w_1,x_1)
      \cdots\gY{i-1}\cdot\\
    &\qquad\qquad \cdot \Y_i(\U(x_1)u_mw_i,x_i)
      \gY {i+1}\cdot\bn
    &\qquad\qquad \cdots \gY n\qL\bn
    &\quad = \sum_{i=1}^n\sum_{m=0}^\infty\sum_{l=1}^\infty
      \frac{(2\pi i)^{m+1}}{m!} \left(x_i\frac{\partial}{\partial x_i}\right)^m
      \left(\frac{x^l}{x_i^l} + \frac{x^{-l}}{x_i^{-l}}\right)\\
    &\qquad\qquad\Tr \Y_1(\U(x_1)w_1,x_1)
      \cdots\gY{i-1}\cdot\bn
    &\qquad\qquad \cdot \Y_i(\U(x_1)u_mw_i,x_i)
      \gY {i+1}\cdot\\
    &\qquad\qquad \cdots \gY n\qL\\
    &\quad\qquad + 2\pi i \sum_{i=1}^n 
      \Y_1(\U(x_1)w_1,x_1)
      \cdots\gY{i-1}\cdot\\
    &\qquad\qquad \cdot \Y_i(\U(x_1)u_0w_i,x_i)
      \gY {i+1} \cdots \gY n\qL
  \end{align*}
  Since the operator $(1-q^{-x{\frac{\partial}{\partial x}}})$ kills constant expressions
  in the variable $x$, the left hand side has no constant term as a series in $x$. This
  implies
  \begin{align*}
    &\sum_{i=1}^n \Tr \gY 1\cdots
      \gY {i-1}\cdot\\
    &\qquad\cdot \Y_i(\U(x_i)u_0w_i,x_i)
      \gY{i+1} \cdots\gY n q^{L(0)}\\
    &\qquad = 0;
  \end{align*}
  therefore,
  \begin{align*}
    &\left(1-q^{-x\frac{\partial}{\partial x}}\right)\Tr Y(\U(x)u,x)
      \gY 1\cdots\gY n \qL\\
    &\quad = \sum_{i=1}^n\sum_{m=0}^\infty\sum_{l=1}^\infty
      \frac{(2\pi i)^{m+1}}{m!} 
      \left((-l)^m\frac{x^l}{x_i^l} + l^m\frac{x^{-l}}{x_i^{-l}}\right)\bn
    &\qquad\qquad\Tr \Y_1(\U(x_1)w_1,x_1)
      \cdots\gY{i-1}\cdot\\
    &\qquad\qquad \cdot \Y_i(\U(x_1)u_mw_i,x_i)
      \gY {i+1}\cdot\bn
    &\qquad\qquad \cdots \gY n\qL.
  \end{align*}
  Then invert the operator $(1-q^{-x\frac{\partial}{\partial x}})$, with the
  appropriate constant term:
  \begin{align*}
    &\Tr Y(\U(x)u,x)\gY 1\cdots\gY n \qL\\
    &\quad = \Tr o(U(1)u)\gY 1\cdots \gY n \qL\bn
    &\quad\qquad + \left(1-q^{-x\frac{\partial}{\partial x}}\right)^{-1}
      \sum_{i=1}^n\sum_{m=0}^\infty\sum_{l=1}^\infty
      \frac{(2\pi i)^{m+1}}{m!} 
      \left((-l)^m\frac{x^l}{x_i^l} + l^m\frac{x^{-l}}{x_i^{-l}}\right)\bn
    &\qquad\qquad\Tr \Y_1(\U(x_1)w_1,x_1)
      \cdots\gY{i-1}\cdot\\
    &\qquad\qquad \cdot \Y_i(\U(x_1)u_mw_i,x_i)
      \gY {i+1}\cdot\bn
    &\qquad\qquad \cdots \gY n\qL\bn
    &\quad = \Tr o(U(1)u)\gY 1\cdots \gY n \qL\\
    &\quad\qquad + \sum_{i=1}^n\sum_{m=0}^\infty\sum_{l=1}^\infty
      \frac{(2\pi i)^{m+1}}{m!} 
      \left(-(-l)^m\frac{q^lx^l}{(1-q^l)x_i^l} + 
      l^m\frac{x^{-l}}{(1-q^l)x_i^{-l}}\right)\bn
    &\qquad\qquad\Tr \Y_1(\U(x_1)w_1,x_1)
      \cdots\gY{i-1}\cdot\\
    &\qquad\qquad \cdot \Y_i(\U(x_1)u_mw_i,x_i)
      \gY {i+1}\cdot\bn
    &\qquad\qquad \cdots \gY n\qL\bn
    &\quad = \Tr o(U(1)u)\gY 1\cdots \gY n \qL\\
    &\quad\qquad + \sum_{i=1}^n\sum_{m=0}^\infty\sum_{l=1}^\infty
      \frac{(2\pi i)^{m+1}}{m!} 
      \left(-(-l)^m\frac{q^l\left(\frac{x_i}{x}\right)^{-l}}{1-q^l} + 
      l^m\frac{\left(\frac{x_i}{x}\right)^{l}}{(1-q^l)}\right)\bn
    &\qquad\qquad\Tr \Y_1(\U(x_1)w_1,x_1)
      \cdots\gY{i-1}\cdot\\
    &\qquad\qquad \cdot \Y_i(\U(x_1)u_mw_i,x_i)
      \gY {i+1}\cdot\bn
    &\qquad\qquad \cdots \gY n\qL\bn
    &\qquad = \Tr o(\U(1)u)\gY 1\cdots \gY n \qL \bn
    &\qquad\qquad + \sum_{i=1}^n\sum_{m\geq 0} P_{m+1}\left(\frac{x_i}{x}; q\right)
    \Tr \Y_1(\U(x_1)w_1,x_1)\cdot\nn
    &\qquad\qquad\qquad\cdots \Y_{i-1}(\U(x_{i-1})w_{i-1},x_{i-1})
      \Y_i(\U(x_i)u_m w_i,x_i)\cdot\nn
    &\qquad\qquad\qquad\cdot \Y_{i+1}(\U(x_{i+1})w_{i+1},x_{i+1})\cdots
      \Y_{n}(\U(x_{n})w_{n},x_{n})q^{L(0)}
  \end{align*}
\epfv

In the rest of this section, we  consider grading-restricted generalized $V$-modules whose weights are real,
and we make the following assumptions which will be proved in the next section:
\begin{itemize}
  \item (Convergence of genus zero correlation functions) For any grading-restricted generalized
    $V$-modules $W_i$, $\tilde W_i$, $i=0, \ldots, n$ with $\tilde W_n = \tilde W_0$,
    for any elements $w_i\in W_i$, $\tilde w_n \in \tilde W_n$,
    $\tilde w_n^\prime \in \tilde W_n^\prime$, and for logarithmic intertwining
    operators $\Y_i\in \binom{\tilde W_{i-1}}{W_i \tilde W_i}$, $i=0\ldots n$,
    the series
    \begin{align*}
      \langle\tilde w_n^\prime, \Y_1(w_1,z_1)\cdots\Y_n(w_n,z_n)\tilde w_n\rangle
    \end{align*}
    is absolutely convergent in the region $|z_1|>|z_2|>\ldots>|z_n|>0$.
  \item (Associativity for $P$-intertwining operators)
    For any grading-restricted generalized $V$-modules $\tilde W_1$, $\tilde W_2$,
    $V$-$P$-bimodules $W_1$, $W_2$, $W_3$, 
    $P$-intertwining operators $\Y_1\in \binom{W_2}{\tilde W_2, W_3}$,
    $\Y_2\in \binom{W_1}{\tilde W_1, W_2}$ there exists
    a $V$-$P$-bimodule $M$, a logarithmic intertwining operator
    $Y_3\in \binom{M}{\tilde W_1, \tilde W_2}$ and a $P$-intertwining operator
    $\Y_4 \in \binom{W_1}{M, W_3}$ such that for any complex numbers
    $z_1, z_2$ with $|z_1|>|z_2|>|z_1-z_2|>0$ and $\tilde w_1 \in \tilde W_1$,
    $\tilde w_2 \in W_2$, $w_3 \in W_3$, $w_1^\prime \in W_1^\prime$ 
    \begin{align*}
      &\langle w_1^\prime, \Y_1(\tilde w_1, z_1)\Y_2(\tilde w_2, z_2)w_3\rangle\\
      &\qquad = \langle w_1^\prime, 
        \Y_4(\Y_3(\tilde w_1, z_1-z_2)\tilde w_2, z_2)w_3\rangle
    \end{align*}
  \item (Commutativity for $P$-intertwining operators) 
    For any grading-restricted generalized $V$-modules $\tilde W_1,$ $\tilde W_2$, $W_1$, $W_2$,
    $W_3$, $P$-intertwining operators $\Y_1\in \binom{W_2}{\tilde W_2, W_3}$,
    $\Y_2\in \binom{W_1}{\tilde W_1, W_2}$there exists
    a $V$-$P$-bimodule $M$ and $P$-intertwining operators
    $Y_3\in \binom{W_1}{\tilde W_2, M}$
    and $\Y_4 \in \binom{M}{\tilde W_1, W_3}$ such that for any
    $\tilde w_1 \in \tilde W_1$,
    $\tilde w_2 \in W_2$, $w_3 \in W_3$, $w_1^\prime \in W_1^\prime$,
    the multivalued analytic function 
    \begin{align*}
      &\langle w_1^\prime, \Y_1(\tilde w_1, z_1)\Y_2(\tilde w_2, z_2)w_3\rangle
    \end{align*}
    in the region $|z_1|>|z_2|>0$ is an analytic continuation of the multivalued
    analytic function
    \begin{align*}
      &\langle w_1^\prime, \Y_3(\tilde w_2, z_2)\Y_4(\tilde w_1, z_1)w_3\rangle
    \end{align*}
    in the region $|z_2|>|z_1|>0$.
\end{itemize}

Our goal is to obtain differential equations for the genus-one correlation
functions. In order to do that, we need to derive formulae for
\begin{align*}
  \Tr \gY 1 \cdots \Y_i(\U(x_i)L(-1)w_i,x_i)\cdots \gY n \qL,
\end{align*}
which is related to the derivative of the formal $q$-trace with respect
to the variable $x_i$; hence, we consider the expression
\begin{align*}
  \Tr \gY 1 \cdots \Y_i(\U(x_i)Y(u, y)w_i,x_i)\cdots \gY n \qL
\end{align*}
Using the $\U(x)$ conjugation property, one can
rewrite this as
\begin{align*}
  &\Tr \gY 1\cdot \nn
  &\qquad \cdots \Y_i(Y(\U(x_ie^{2\pi i y}u, x_i(e^{2 \pi i y}-1))
    \U(x_i)w_i,x_i)\cdot\nn
  &\qquad\qquad \cdots \gY n \qL
\end{align*}
then, in order to use (\ref{identity1}), one rewrites the iterate as a product using
associativity for intertwining operators.

\begin{lemma}
  Let $\Y$ be a logarithmic intertwining operator of type $\binom{W_0}{W W_1}$
  for grading-restricted generalized $V$-modules $W, W_0, W_1$. Then for any
  $w_0^\prime \in W_0^\prime, w_1\in W_1, w\in W$, and for any complex number
  $z$ satisfying $|q_z| > 1 > |q_z - 1| > 0$,
  \begin{align} \label{rf1}
    &\langle w_0^\prime, \Y(Y(\U(xq_z)u, x(q_z-1))\U(x)w, x) w_1\rangle \nn
    &\qquad = \langle w_0^\prime, Y(\U(xq_z)u, xq_z)\Y(\U(x)w, x)w_1\rangle
  \end{align}
\end{lemma}
\proof
  Using associativity for intertwining operators, we see that
  \begin{align}\label{rf2}
    &\langle w_0^\prime, \Y(Y(\U(z_1q_z)u, z_1(q_z-1))\U(z_1)w, z_1) w_1\rangle \nn
    &\qquad = \langle w_0^\prime, Y(\U(z_1q_z)u, z_1q_z)\Y(\U(z_1)w, z_1)w_1\rangle
  \end{align}
  holds for any complex numbers $z, z_1$ in the region
  $|z_1q_z| > |z_1| > |z_1(q_z - 1)| > 0$, or whenever $z_1\neq 0$ and
  $|q_z| > 1 > |q_z - 1| > 0$.
  Then for a fixed $z$ in the above region, we have two series in powers of $z_1$
  (not necessarily integral), and $\log z_1$; by our assumption on the modules,
  these powers must form  a unique expansion set, and thus the coefficients of the
  two series must be equal. Then, we can replace the complex variable $z_1$ in
  (\ref{rf2}) with the formal variable $x$, which concludes the proof.
\epfv
\begin{prop}
  Consider grading-restricted generalized $V$-modules $W_i$ and $V$-$P$-bimodules
  $\tilde W_i$ for $i=1,\ldots,n$, with $\tilde W_0 = \tilde W_n$, and 
  $P$-intertwining operators $\Y_i$ of type $\binom{\tilde W_{i-1}}{W_i \tilde W_i}$
  for $i=1,\ldots,n$. Moreover, suppose $\tilde W_0$ is projective as right $P$-module.
  Then for any $v\in V$, $w_i\in W_i$, and any integer $j$, $1\leq j \leq n$,
  \begin{align} \label{Y-iterate-1}
    &\Tr \gY 1\cdots \gY{j-1}\Y_j(\U(x_j)Y(v,y)w_j,x_j)\bn
    &\qquad\qquad \cdot\gY{j+1} \cdots \gY n \qL\bn
    &\quad = \sum_{m\geq 0}(-1)^{m+1}
      \left( \tilde \wp_{m+1}(-y;q)
      +\frac{\partial^m}{\partial y^m}(\tilde G_2(q)y + \pi i)\right)\cdot\bn
    &\qquad\qquad \Tr \gY 1 \cdots \gY{j-1}\cdot\bn
    &\qquad\qquad \cdot \Y_j(\U(x_j)v_mw_j,x_j)\gY{j+1}\cdot\bn
    &\qquad\qquad \cdots \gY{n} \qL\bn
    &\quad\quad + \sum_{i\neq j}\sum_{m\geq 0} P_{m+1}\left(
      \frac{x_i}{x_je^{2\pi i y}}; q \right)\cdot\bn
    &\qquad\qquad \cdot\Tr \gY{1}\cdots \gY{i-1}\cdot\bn
    &\qquad\qquad \cdot \Y_i(\U(x_i)v_mw_i, x_i)\gY{i+1}
      \cdots \gY{n} \qL\bn
    &\quad\quad + \Tr o(\U(1)v)\gY{1} \cdots \gY{n}\qL.
  \end{align}
\end{prop}
\proof
  By induction on $j$; when $j=1$, by Lemma \ref{U-conj},
  \begin{align*}
    &\Tr \Y_1(\U(x_1)Y(v,y)w_1,x_1)\gY 2 \cdots \gY n \qL\bn
    &\quad = \Tr \Y_1(Y(\U(x_1e^{2\pi i y})v,x_1(e^{2\pi i y}-1))\U(x_1)w_1,x_1)\cdot\bn
    &\qquad\qquad \cdot \gY 2 \cdots \gY n \qL
  \end{align*}
  and using (\ref{rf1}), for any complex number $z$ such that
  $|q_z| > 1 > |q_z-1| > 0$,
  \begin{align*}
    &\Tr \Y_1(\U(x_1)Y(v,z)w_1,x_1)\gY 2 \cdots \gY n \qL\bn
    &\quad = \Tr \Y_1(Y(\U(x_1q_z)v,x_1(q_z-1))\U(x_1)w_1,x_1)\cdot\bn
    &\qquad\qquad \cdot \gY 2 \cdots \gY n \qL\bn
    &\quad = \Tr Y(\U(x_1q_z)v,x_1q_z)\Y_1(\U(x_1)w_1,x_1)\cdot\bn
    &\qquad\qquad \cdot \gY 2 \cdots \gY n \qL
  \end{align*}
  Now by (\ref{identity1}) with $x = x_1q_z$, we get
  \begin{align*}
    &\Tr \Y_1(\U(x_1)Y(v,z)w_1,x_1)\gY 2 \cdots \gY n \qL\bn
    &\quad = \sum_{m\geq 0} P_{m+1} \left(\frac{1}{q_z}; q\right)
      \Tr \Y_1(\U(x_1)v_mw_1,x_1)\cdot \bn
    &\qquad\qquad \cdot \gY 2 \cdots \gY n\qL \bn
    &\quad\quad + \sum_{i=2}^m \sum_{m\geq 0} P_{m+1}
      \left(\frac{x_i}{x_1q_z}; q\right) \Tr\gY 1\cdots \gY{i-1}\cdot\bn
    &\qquad\qquad \cdot \Y_i(\U(x_1)v_mw_i,x_i)\cdot \gY {i+1} \cdots
      \gY n \qL\bn
    &\qquad\quad + \Tr o(\U(1)v)\gY 1 \cdots \gY n \qL\bn
    &\quad = \sum_{m\geq 0} (-1)^{m+1}
      \left(\tilde\wp_{m+1}(-z;q) + \frac{\partial^m}{\partial z^m}
      (\tilde G_2(q)z + \pi i)\right)
      \Tr \Y_1(\U(x_1)v_mw_1,x_1)\cdot \bn
    &\qquad\qquad \cdot \gY 2 \cdots \gY n\qL \bn
    &\quad\quad + \sum_{i=2}^m \sum_{m\geq 0} P_{m+1}
      \left(\frac{x_i}{x_1q_z}; q\right) \Tr\gY 1\cdots \gY{i-1}\cdot\bn
    &\qquad\qquad \cdot \Y_i(\U(x_1)v_mw_i,x_i)\cdot \gY {i+1} \cdots
      \gY n \qL\bn
    &\qquad\quad + \Tr o(\U(1)v)\gY 1 \cdots \gY n \qL
  \end{align*}
  which proves the base case. Now suppose the result holds for $j\geq 1$,
  and use commutativity for intertwining operators.
\epfv

\begin{prop} Using the same notation as in the previous proposition, for all
  $v\in V$ and $l\geq 1$, we have
  \begin{align} \label{rf3}
    &\Tr \gYq{1} \cdots \gYq{j-1}
      \Y_j(\U(q_{z_j})v_{-l}w_j, q_{z_j})\cdot \bn
    &\qquad\qquad \cdot \gYq{i+1}\cdots \gYq{n} \qL \bn
    &\qquad = \sum_{k=1}^\infty (-1)^{l+1} \binom{2k+1}{l-1}
      \tilde G_{2k+2}(q) \Tr \gYq{1} \cdot \bn
    &\qquad\qquad\cdots\gYq{j-1}\Y_j(\U(q_{z_j})v_{2k+2-l}w_j, q_{z_j})\cdot\bn
    &\qquad\qquad\cdot \gYq{i+1}\cdots \gYq{n}\qL\bn
    &\qquad\quad + \sum_{i\neq j}\sum_{m=0}^\infty(-1)^{m+l}\binom{-m-1}{l-1}
      \tilde \wp_{m+l}(z_i-z_j;q) \Tr \gYq{1} \cdot \bn
    &\qquad\qquad\cdots\gYq{i-1}\Y_i(\U(q_{z_i})v_{m}w_i, q_{z_i})\cdot\bn
    &\qquad\qquad\cdot \gYq{i+1}\cdots \gYq{n}\qL\bn
    &\qquad\quad + \delta_{l,1}\tilde G_2(q)\sum_{i=1}^n
      \Tr \gYq{1}\cdots\gYq{i-1}\cdot\bn
    &\qquad\qquad\cdot \Y_i(\U(q_{z_i})(v_1+v_0z_i)w_i, q_{z_i})\gYq{i+1}\cdot\bn
    &\qquad\qquad\cdots \gYq{n}\qL\bn
    &\qquad\quad + \delta_{l,1} \Tr o(\U(1)v)\gYq{1}\cdots\gYq{n}\qL
  \end{align}
\end{prop}
\proof
Notice that the coefficients of (\ref{Y-iterate-1}) as a series in the formal variables
$q$, $\log q$ are absolutely convergent in the region $|q_{z_1}| > \ldots > |q_{z_n}| > 0$.
Since the $q$ coefficients of $\tilde \wp_m(z, q)$ are absolutely convergent
when $|z| < 1$, we can substitute $y = z$, $x_i = q_{z_i}$, for $i=1,\ldots, n$ in
(\ref{Y-iterate-1}),
\begin{align*} 
  &\Tr \gYq 1\cdots \gYq{j-1}\Y_j(\U(q_{z_j})Y(v,z)w_j,q_{z_j})\bn
    &\qquad\qquad \cdot\gYq{j+1} \cdots \gYq n \qL\bn
    &\quad = \sum_{m\geq 0}(-1)^{m+1}
      \left( \tilde \wp_{m+1}(-z;q)
      +\frac{\partial^m}{\partial z^m}(\tilde G_2(q)z + \pi i)\right)\cdot\bn
    &\qquad\qquad \Tr \gYq 1 \cdots \gYq{j-1}\cdot\bn
      &\qquad\qquad \cdot \Y_j(\U(q_{z_j})v_mw_j,q_{z_j})\gYq{j+1}\cdot\bn
    &\qquad\qquad \cdots \gYq{n} \qL\bn
    &\quad\quad + \sum_{i\neq j}\sum_{m\geq 0} 
      \tilde (-1)^{m+1}\Bigg(\wp_{m+1}(z_i - z_j - z; q)\\
    &\hspace{6cm} + (-1)^{m+1}\der{m}{z_i}(\tilde G_2(q)(z_i - z_j - z) - \pi i)\Bigg)
      \cdot\bn
    &\qquad\qquad \cdot\Tr \gYq{1}\cdots \gYq{i-1}\cdot\bn
      &\qquad\qquad \cdot \Y_i(\U(q_{z_i})v_mw_i, q_{z_i})\gYq{i+1}
      \cdots \gYq{n} \qL\bn
    &\quad\quad + \Tr o(\U(1)v)\gYq{1} \cdots \gYq{n}\qL.
  \end{align*}
  The result follows by taking the coefficient of $z^{l-1}$ and using
  the $q$ expansion of $\tilde \wp(z; q)$ and (\ref{identity2}).
\epfv

\noindent
Taking $v = \omega$ and $l=1$ in (\ref{rf3}), one obtains
\begin{align*}
    &\Tr \gYq{1} \cdots \gYq{j-1}
      \Y_j(\U(q_{z_j})\omega_{-1}w_j, q_{z_j})\cdot \bn
    &\qquad\qquad \cdot \gYq{i+1}\cdots \gYq{n} \qL \bn
    &\qquad = \sum_{k=1}^\infty 
      \tilde G_{2k+2}(q) \Tr \gYq{1} \cdot \bn
    &\qquad\qquad\cdots\gYq{j-1}\Y_j(\U(q_{z_j})\omega_{2k+1}w_j, q_{z_j})\cdot\bn
    &\qquad\qquad\cdot \gYq{i+1}\cdots \gYq{n}\qL\bn
    &\qquad\quad + \sum_{i\neq j}\sum_{m=0}^\infty(-1)^{m+1}
      \tilde \wp_{m+1}(z_i-z_j;q) \Tr \gYq{1} \cdot \bn
    &\qquad\qquad\cdots\gYq{i-1}\Y_i(\U(q_{z_i})\omega_{m}w_i, q_{z_i})\cdot\bn
    &\qquad\qquad\cdot \gYq{i+1}\cdots \gYq{n}\qL\bn
    &\qquad\quad + \tilde G_2(q)\sum_{i=1}^n
      \Tr \gYq{1}\cdots\gYq{i-1}\cdot\bn
    &\qquad\qquad\cdot \Y_i(\U(q_{z_i})(\omega_1+\omega_0z_i)w_i, q_{z_i})\gYq{i+1}\cdot\bn
    &\qquad\qquad\cdots \gYq{n}\qL\bn
    &\qquad\quad + \Tr o(\U(1)\omega)\gYq{1}\cdots\gYq{n}\qL
\end{align*}
now since $w_k = L(k-1)$, and since 
$\U(1)\omega = (2\pi i)^2\left(\omega - \frac{c}{24}\one\right)$, we obtain
\begin{align} \label{rf4}
    &\Tr \gYq{1} \cdots \gYq{j-1}
  \Y_j(\U(q_{z_j})L(-2)w_j, q_{z_j})\cdot \bn
    &\qquad\qquad \cdot \gYq{i+1}\cdots \gYq{n} \qL \bn
    &\qquad = \sum_{k=1}^\infty 
      \tilde G_{2k+2}(q) \Tr \gYq{1} \cdot \bn
    &\qquad\qquad\cdots\gYq{j-1}\Y_j(\U(q_{z_j})L(2k)w_j, q_{z_j})\cdot\bn
    &\qquad\qquad\cdot \gYq{i+1}\cdots \gYq{n}\qL\bn
    &\qquad\quad + \sum_{i\neq j}\sum_{m=0}^\infty(-1)^{m+1}
      \tilde \wp_{m+1}(z_i-z_j;q) \Tr \gYq{1} \cdot \bn
    &\qquad\qquad\cdots\gYq{i-1}\Y_i(\U(q_{z_i})L(m-1)w_i, q_{z_i})\cdot\bn
    &\qquad\qquad\cdot \gYq{i+1}\cdots \gYq{n}\qL\bn
    &\qquad\quad + \tilde G_2(q)\sum_{i=1}^n
      \Tr \gYq{1}\cdots\gYq{i-1}\cdot\bn
    &\qquad\qquad\cdot \Y_i(\U(q_{z_i})(L(0)+L(-1)z_i)w_i, q_{z_i})\gYq{i+1}\cdot\bn
    &\qquad\qquad\cdots \gYq{n}\qL\bn
    &\qquad\quad + (2\pi i)^2
      \Tr \left(L(0) - \frac{c}{24}\right)\cdot\bn
    &\qquad\qquad \cdot\gYq{1}\cdots\gYq{n}\qL.
\end{align}

\subsection{Duality properties for $P$-intertwining operators}
\label{s-duality}

In this section we derive associativity and commutativity properties
for $P$-intertwining operators, under the same assumptions
as Proposition \ref{HLZassociativity}: we show that in the statement
of these two properties, all the modules can be taken to
be $V$-$P$-bimodules and all intertwining operators to be
$P$-intertwining operators. In particular, under those assumptions
all identities in the previous sections hold.

Let $W_1$ be a grading-restricted generalized $V$-module and $W_2$ be a $V$-$P$-bimodule. We can
define a right action of $P$ on $W_1\boxtimes_{P(z)} W_2$ the following way:
for any $p\in P$, consider the map
\begin{align*}
    \text{Id}_{W_1}\otimes p : W_1\otimes W_2 &\rightarrow W_1\otimes W_2\\
    w_1\otimes w_2 &\mapsto w_1 \otimes (w_2p);
\end{align*}
then, by Proposition \ref{tensor-morphism}, there exists a unique $V$-homomorphism
$\text{Id}_{W_1}\boxtimes_{P(z)} p$ of $W_1\boxtimes_{P(z)} W_2$
such that
\begin{equation}\label{sr1}
  \boxtimes_{P(z)}\circ (\text{Id}_{W_1}\otimes p) =
  \overline{(\text{Id}_{W_1}\boxtimes_{P(z)} p)} \circ \boxtimes_{P(z)}.
\end{equation}
so we let $p$ act on $W_1\boxtimes_{P(z)} W_2$ by
\begin{align*}
  wp = (\text{Id}_{W_1}\boxtimes_{P(z)}p)(w)
\end{align*}
for any $w\in W_1\boxtimes_{P(z)} W_2$, and extend it to the formal completion
$\overline {W_1\boxtimes_{P(z)} W_2}$; thus (\ref{sr1}) becomes
\begin{equation}
  \boxtimes_{P(z)}(w_1\otimes (w_2p)) =
  \boxtimes_{P(z)}(w_1\otimes w_2)p
\end{equation}
for all $w_1\in W_1$ and $w_2 \in W_2$.
Thus, $W_1\boxtimes_{P(z)}W_2$ can be seen naturally
as a $V$-$P$-bimodule; moreover, the intertwining
map $\boxtimes_{P(z)}$ commutes with the action of $P$.
In particular, $P$ commutes with the intertwining operator $\Y_{\boxtimes_{P(z)},0}$.

\begin{prop}
  Let $\Y_1 \in \binom{W_4}{W_1 M}$ and $\Y_2\in \binom{M}{W_2W_3}$ be two
  logarithmic intertwining operators, where $W_1, W_2, W_3$ are $V$-$P$-bimodules,
  $M$ is a grading-restricted generalized $V$-module, and $P$ commutes with $\Y_1$ and $\Y_2$.
  Let $\Y^1 \in \binom{W_4}{WW_3}$ be the logarithmic intertwining operator
  as in Proposition \ref{HLZassociativity} (1.) (here $W= W_1\boxtimes_{P(z_0)}W_2$).
  Then $P$ commutes with $\Y^1$.
\end{prop}
\proof
  Let $p$ be any element in $P$; for any two $V$-$P$-bimodules $W^2, W^3$,
 grading-restricted  generalized $V$-module $W^1$, and logarithmic intertwining operator
  $\Y\in \binom{W^3}{W^1 W^2}$, consider the logarithmic intertwining
  operators $\Y_p, \Y^p \in \binom{W^3}{W^1W^2}$ defined by
  \begin{eqnarray*}
    \Y_p(w^{(1)}, x)w^{(2)} &= \Y(w^{(1)}, x)(w^{(2)}p)\\
    \Y^p(w^{(1)}, x)w^{(2)} &= (\Y(w^{(1)}, x)w^{(2)})p
  \end{eqnarray*}
  for all $w^{(1)}\in W^1, w^{(2)} \in W^2$.
  It is clear that $\Y_p$ (resp. $\Y^p$) is a logarithmic intertwining operator
  since $p$ acts on $W^2$ (resp. $W^3$) as a $V$-module homomorphism.
  Then $P$ commutes with $\Y$ if and only if $\Y_p = \Y^p$ for all
  $p\in P$.
  Consider $(\Y^1)_p$ and $(\Y^1)^p$: then for $w_{(1)}\in W_1$,
  $w_{(2)}\in W_2$,
  $w_{(3)}\in W_3$,
  $w_{(4)}^\prime\in W_4^\prime$, and for all $z_0, z_1, z_2$ such that
  $z_0 = z_1 - z_2$, $|z_1|>|z_2|>|z_0|>0$,
  \begin{align*}
    &\langle w_{(4)}^\prime, (\Y^1)^p
      (\Y_{\boxtimes_{P(z_0)},0}(w_{(1)},z_0)w_{(2)},z_2)w_{(3)}\rangle\\
    &\qquad = \langle w_{(4)}^\prime, (\Y^1
      (\Y_{\boxtimes_{P(z_0)},0}(w_{(1)},z_0)w_{(2)},z_2)w_{(3)})p\rangle\bn
    &\qquad = \langle pw_{(4)}^\prime, \Y^1
      (\Y_{\boxtimes_{P(z_0)},0}(w_{(1)},z_0)w_{(2)},z_2)w_{(3)}\rangle\\
    &\qquad = \langle pw_{(4)}^\prime, \Y_1(w_{(1)}, z_1)
      \Y_2(w_{(2)},z_2)w_{(3)}\rangle\\
    &\qquad = \langle w_{(4)}^\prime, (\Y_1(w_{(1)}, z_1)
      \Y_2(w_{(2)},z_2)w_{(3)})p\rangle\bn
    &\qquad = \langle w_{(4)}^\prime, \Y_1(w_{(1)}, z_1)
      (\Y_2(w_{(2)},z_2)w_{(3)})p\rangle\\
    &\qquad = \langle w_{(4)}^\prime, \Y_1(w_{(1)}, z_1)
      (\Y_2)^p(w_{(2)},z_2)w_{(3)}\rangle.
  \end{align*}
  Similarly,
  \begin{align*}
    &\langle w_{(4)}^\prime, (\Y^1)_p
      (\Y_{\boxtimes_{P(z_0)},0}(w_{(1)},z_0)w_{(2)},z_2)w_{(3)}\rangle\\
    &\qquad = \langle w_{(4)}^\prime, \Y^1(
      \Y_{\boxtimes_{P(z_0)},0}(w_{(1)},z_0)w_{(2)},z_2)(w_{(3)}p)\rangle\\
    &\qquad = \langle w_{(4)}^\prime, \Y_1(w_{(1)}, z_1)
      \Y_2(w_{(2)},z_2)(w_{(3)}p)\rangle\\
    &\qquad = \langle w_{(4)}^\prime, \Y_1(w_{(1)}, z_1)
      (\Y_2)_p(w_{(2)},z_2)w_{(3)}\rangle\\
    &\qquad = \langle w_{(4)}^\prime, \Y_1(w_{(1)}, z_1)
      (\Y_2)^p(w_{(2)},z_2)w_{(3)}\rangle.
  \end{align*}
  Thus, by uniqueness in Proposition \ref{HLZassociativity} (2.) applied to the
 logarithmic  intertwining operators $\Y_1$ and $(\Y_2)^p$ we see that
  $(\Y^1)^p = (\Y^1)_p$ and therefore $P$ commutes with $\Y^1$.
\epfv
\begin{prop}
  Using the notation of Proposition \ref{HLZassociativity} (2.), suppose the modules
  $W_3$ and $W_4$ are $V$-$P$-bimodules, and the logarithmic intertwining operator $\Y^1$
  commutes with $P$. Then the logarithmic intertwining operator $\Y_1$ also
  commutes with $P$.
\end{prop}
\proof
  Note that since $W_3$ is a $V$-$P$-bimodule, the right action of $P$ on $W_3$
  defines a right action of $P$ on $W_2\boxtimes_{P(z)} W_3$, and $\Y_{\boxtimes_{P(z_2)},0}$
  commutes with $P$.
  Now for $w_{(1)}\in W_1$, $w_{(2)}\in W_2$, $w_{(3)}\in W_3$ and 
  $w_{(4)}^\prime \in W_4^\prime$, and for complex numbers $z_0, z_1, z_2$ such
  that $|z_1| > |z_2| > |z_0| > 0$, $z_0 = z_1 - z_2$, for any $p\in P$
  \begin{align*}
    &\langle w_{(4)}^\prime, (\Y_1)^p(w_{(1)}, z_1)
      \Y_{\boxtimes_{P(z_2)},0}(w_{(2)},z_2)w_{(3)}\rangle\bn
    &\qquad = \langle w_{(4)}^\prime, (\Y_1(w_{(1)}, z_1)
      \Y_{\boxtimes_{P(z_2)},0}(w_{(2)},z_2)w_{(3)})p\rangle\bn
    &\qquad = \langle pw_{(4)}^\prime, \Y_1(w_{(1)}, z_1)
      \Y_{\boxtimes_{P(z_2)},0}(w_{(2)},z_2)w_{(3)}\rangle\bn
    &\qquad = \langle pw_{(4)}^\prime, 
      \Y^1(\Y^2(w_{(1)},z_0)w_{(2)},z_2)w_{(3)}\rangle\bn
    &\qquad = \langle w_{(4)}^\prime, 
      (\Y^1(\Y^2(w_{(1)},z_0)w_{(2)},z_2)w_{(3)})p\rangle\bn
    &\qquad = \langle w_{(4)}^\prime, 
      (\Y^1)^p(\Y^2(w_{(1)},z_0)w_{(2)},z_2)w_{(3)}\rangle.
  \end{align*}
  Similarly, since $\Y_{\boxtimes_{P(z_2)},0}$ commutes with $P$,
  \begin{align*}
    &\langle w_{(4)}^\prime, (\Y_1)_p(w_{(1)}, z_1)
      \Y_{\boxtimes_{P(z_2)},0}(w_{(2)},z_2)w_{(3)}\rangle\bn
    &\qquad = \langle w_{(4)}^\prime, \Y_1(w_{(1)}, z_1)
      (\Y_{\boxtimes_{P(z_2)},0}(w_{(2)},z_2)w_{(3)})p\rangle\bn
    &\qquad = \langle w_{(4)}^\prime, \Y_1(w_{(1)}, z_1)
      \Y_{\boxtimes_{P(z_2)},0}(w_{(2)},z_2)(w_{(3)}p)\rangle\bn
    &\qquad = \langle w_{(4)}^\prime, 
      \Y^1(\Y^2(w_{(1)},z_0)w_{(2)},z_2)(w_{(3)}p)\rangle\bn
    &\qquad = \langle w_{(4)}^\prime, 
      (\Y^1)_p(\Y^2(w_{(1)},z_0)w_{(2)},z_2)w_{(3)}\rangle\bn
    &\qquad = \langle w_{(4)}^\prime, 
      (\Y^1)^p(\Y^2(w_{(1)},z_0)w_{(2)},z_2)w_{(3)}\rangle.
  \end{align*}
  The conclusion thus follows from uniqueness in Proposition \ref{HLZassociativity}
  part (2.) applied to the logarithmic intertwining operators $(\Y^1)^p$ and $\Y^2$.
\epfv

We summarize these results in the following.
\begin{thm}[Associativity for $P$-intertwining operators] \label{assoc1} 
  Let $W_1, W_2, W_3$ be $V$-$P$-bimodules.

  (i) Let $M$ be a $V$-$P$-bimodule and
  $\Y_1$, $\Y_2$ be two $P$-intertwining operators of types
  $\binom{W_4}{W_1 M}$, $\binom{M}{W_2W_3}$, respectively. Then there exist a
  grading-restricted generalized module $W$, a $P$-intertwining operator $\Y^1$ of type $\binom{W_4}{WW_3}$
  and a logarithmic intertwining operator $\Y^2$ of type $\binom{W}{W_1 W_2}$ such that
  for all $w_{(1)}\in W_1$, $w_{(2)}\in W_2$, $w_{(3)}\in W_3$, $w_{(4)}\in W_4^\prime$,
  \begin{align*} 
    &\langle w_{(4)}^\prime, \Y_1(w_{(1)}, z_1)\Y_2(w_{(2)}, z_2) w_{(3)}\rangle\\
    &\qquad = \langle w_{(4)}^\prime, \Y^1(\Y^2(w_{(1)}, z_1-z_2)w_{(2)}, z_2) w_{(3)}\rangle
  \end{align*}
  for all $z_1, z_2$ such that $|z_1| > |z_2| > |z_1 - z_2|$.

  (ii) Let $W$ be a grading-restricted generalized $V$-module, $\Y^1$ a logarithmic intertwining
  operator of type $\binom{W_4}{WW_3}$, and $\Y^2$ a $P$-intertwining operator of type
  $\binom{W}{W_1 W_2}$. Then there exists a $V$-$P$-bimodule $M$ and $P$-intertwining
  operators $\Y_1$, $\Y_2$ of types $\binom{W_4}{W_1 M}$, $\binom{M}{W_2W_3}$, respectively
  such that the same conclusion as in (i) holds.
\end{thm}

We now state the commutativity property.

\begin{thm}[Commutativity for $P$-intertwining operators] \label{commut1}
  Let $\Y_1$, $\Y_2$ be $P$-intertwining operators of types $\binom{W_4}{W_1 M}$,
  $\binom{M}{W_2 W_3}$, respectively, for grading-restricted generalized $V$-modules
  $W_1, W_2$ and $V$-$P$-bimodules $W_3, W_4, M$. 
  Then there exist a $V$-$P$-bimodule $M_1$ and logarithmic intertwining
  operators $\Y_3, \Y_4$ of type $\binom{W_4}{W_2 M_1}$, $\binom{M_1}{W_1 W_3}$
  respectively, commuting with the action of $P$, such that for any
  $w_{(1)}\in W_1$, $w_{(2)}\in W_2$, $w_{(3)}\in W_3$,
  $w_{(4)}^\prime \in W_4^\prime$, the multivalued analytic function
  \begin{equation*}
    \langle w_{(4)}^\prime, \Y_1(w_{(1)}, z_1)\Y_2(w_{(2)}, z_2) w_{(3)}\rangle
  \end{equation*}
  on the region $|z_1|>|z_2|>0$ and the multivalued analytic function
  \begin{equation*}
    \langle w_{(4)}^\prime, \Y_3(w_{(2)}, z_2)\Y_4(w_{(1)}, z_1) w_{(3)}\rangle
  \end{equation*}
  on the region $|z_2|>|z_1|>0$ are analytic extensions of each other.
\end{thm}
\proof
  Using associativity twice, by the previous propositions one sees that all the outer
  intertwining operators commute with $P$.
  Recall (\cite{HLZ2}) the operator
  $\Omega_r : \mathcal V_{W^1W^2}^{W^3}
  \rightarrow \mathcal V_{W^2W^1}^{W^3}$ defined by
  \begin{align*}
    \Omega_r(\Y)(w_{(2)}, x)w_{(1)} = 
    e^{xL(-1)} \Y(w_{(1)}, e^{(2r+1)\pi i} x) w_{(2)}
  \end{align*}
  for any $r\in \Z$;
  then $\Omega_{-r-1}(\Omega_r(\Y)) = 
  \Omega_{r}(\Omega_{-r-1}(\Y)) = \Y$ for any $\Y \in \mathcal V_{W^1W^2}^{W^3}$.
  
  We know by Theorem \ref{assoc1} that there exist a module $M$,
  an intertwining operator $\Y^1$ and a $P$-intertwining operator $\Y^2$ such that
  \begin{align*}
    &\langle w_{(4)}^\prime, \Y_1(w_{(1)}, z_1)
      \Y_2(w_{(2)}, z_2) w_{(3)}\rangle\\
    &\qquad = \langle w_{(4)}^\prime, \Y^2(\Y^1(w_{(1)}, z_0)
      w_{(2)}, z_2) w_{(3)}\rangle;
  \end{align*}
  now substituting
  \begin{align*}
    \Y^1(w_{(1)}, x_0)w_{(2)}
      &= \Omega_0(\Omega_{-1}(\Y^1))(w_{(1)}, x_0)w_{(2)}\\
      &= e^{x_0L(-1)}\Omega_{-1}(\Y^1)(w_{(2)}, e^{\pi i}x_0)w_{(1)}
  \end{align*}
  we obtain (in the region $|z_1| > |z_2| > |z_1 - z_2|$)
  \begin{align*}
    &\langle w_{(4)}^\prime, \Y_1(w_{(1)}, z_1)
      \Y_2(w_{(2)}, z_2) w_{(3)}\rangle\\
    &\qquad = \langle w_{(4)}^\prime, 
      \Y^2(\Omega_1(\Y^1)(w_{(2)},
      e^{\pi i}z_0) w_{(1)}, z_2 + z_0) w_{(3)}\rangle;
  \end{align*}
  which is an extension of $\langle w_{(4)}^\prime, 
      \Y^2(\Omega_1(\Y^1)(w_{(2)},
      z_2 - z_1) w_{(1)}, z_1) w_{(3)}\rangle$
  defined on the region $|z_1| > |z_2 - z_1| > 0$. Now by Theorem \ref{assoc1} (ii),
  we know that there exist a $V$-$P$-bimodule $M_1$ and $P$ intertwining operators
  $\Y_3$, $\Y_4$ of type $\binom{W_4}{W_2 M_1}$, $\binom{M_1}{W_1 W_3}$
  respectively, such that
  \begin{align*}
    &\langle w_{(4)}^\prime, 
      \Y^2(\Omega_1(\Y^1)(w_{(2)},
      z_2 - z_1) w_{(1)}, z_1) w_{(3)}\rangle\\
    &\qquad = \langle w_{(4)}^\prime, \Y_3(w_{(2)}, z_2)\Y_4(w_{(1)}, z_1) w_{(3)}\rangle
  \end{align*}
  This concludes the proof.
\epfv

\section{Genus-one correlation functions} \label{sec:2}

\subsection{Differential Equations}
\label{s-diffeq}

In this section we derive a system of differential equations satisfied by
the formal $q$-traces, using the identities obtained in the previous
section. The main technical assumption used in this section is
the $C_2$-cofiniteness of  $V$, which is used
to prove that certain particular modules over a ring of functions are finitely
generated.

In this section, we assume that the vertex operator algebra $V$ is $C_2$-cofinite.
Then by Lemma 2.11 and Proposition 4.3 in \cite{H6}, every grading-restricted 
generalized $V$-module is $C_2$-cofinite and is $\mathbb{R}$-graded.
We denote by $G$ the space of all multivalued analytic functions defined
on the region $|z_1|>|z_2|>\ldots > |z_n| > 0$ with preferred branches on
$|z_1|>|z_2|>\ldots > |z_n| > 0$, $0\leq \arg z_i < 2\pi$ for
$i=1\ldots n$. For any such function $f(z_1,\ldots,z_n)$, the function
$f(q_{z_1}, \ldots, q_{z_n})$ is a multivalued analytic function
defined when $|q_{z_1}| > \ldots > |q_{z_n}|$.
We denote by $G_q$ the space of such functions.
Let 
$$R = \C[\tilde G_4(q), \tilde G_6(q), \tilde \wp_2(z_i-z_j;q), \tilde \wp_3(z_i-z_j;q)].$$
For $V$-modules $W_1, W_1, \ldots, W_n$ we denote by $T$ the graded $R$-module
$$
  R\otimes W_1 \otimes \ldots \otimes W_n,
$$
with grading induced by the grading by generalized eigenvalues on $W_1,\ldots, W_n$.
Denote by $T_r$ the homogeneous subspace of degree $r\in \R$; moreover, define a
filtration on $T$ by $F(T)_r = \coprod_{s\leq r} T_s$.
Let $J$ be the $R$-submodule of $T$ generated by the elements
\begin{align*}
  &A_j(v; w_1, \ldots, w_n)\\
  &\qquad = 1\otimes w_1\otimes\ldots\otimes w_{j-1}\otimes v_{-2}w_j\otimes
    w_{j+1}\otimes\ldots\otimes w_n\\
  &\qquad + \sum_{k=1}^\infty (2k+1)\tilde G_{2k+2}(q)\otimes
    w_1\otimes\ldots\otimes w_{j-1}\otimes v_{2k}w_j\otimes
    w_{j+1}\otimes\ldots\otimes w_n\\
  &\qquad + \sum_{i\neq j}\sum_{m=0}^\infty (-1)^m(m+1)
    \tilde \wp_{m+2}(z_i - z_j)\otimes\\
  &\qquad \qquad w_1\otimes\ldots\otimes w_{j-1}\otimes v_{m}w_j\otimes
    w_{j+1}\otimes\ldots\otimes w_n
\end{align*}
for $v\in V$, $w_i\in W_i$, $i=1\ldots,n$, and $1\leq j \leq n$. The filtration
on $T$ induces one on $J$, and we will denote by $F(J)_r$ the $r$-th subspace in
the filtration.
\begin{prop}\label{pp}
  There exists $N\in \R$ such that for any $r\in \R$, $F(T)_r = F(T)_N + F(J)_r$.
\end{prop}
\proof
  Let $N$ such that
  \begin{equation*}
    \coprod_{r>N} T_r \subseteq \sum_{i=1}^n
      R\otimes W_1\otimes\ldots
        W_{i-1}\otimes C_2(W_i)\otimes W_{i+1}\otimes\ldots W_n;
  \end{equation*}
  then clearly if $r\leq N$, $F(T)_r \subseteq F(T)_N = F(T)_N + F(J)_r$. We now
  prove on induction on $k\in \N$ that if $r = N+k$, then
  $F(T)_r = F(T)_n + F(J)_r$.
  Let $r = N + k + 1$ and let $t \in F(T)_r$. By definition of $N$, we can assume
  $$
    t = 1\otimes w_1\otimes\cdots \otimes
    w_{i-1}\otimes v_{-2}w_i\otimes w_{i+1}\otimes \ldots \otimes w_n
  $$
  for some $v\in V$ and $w_j \in W_j$, $j=1\ldots, n$. Observe that
  the element $A(v; w_1,\ldots, w_n)$ belongs to $F(J)_r$ and
  \begin{align*}
    S = &\sum_{k=1}^\infty (2k+1)\tilde G_{2k+2}(q)\otimes
      w_1\otimes\ldots\otimes w_{j-1}\otimes v_{2k}w_j\otimes
      w_{j+1}\otimes\ldots\otimes w_n\\
    &\qquad + \sum_{i\neq j}\sum_{m=0}^\infty (-1)^m(m+1)
      \tilde \wp_{m+2}(z_i - z_j)\otimes\\
    &\qquad \qquad w_1\otimes\ldots\otimes w_{j-1}\otimes v_{m}w_j\otimes
      w_{j+1}\otimes\ldots\otimes w_n
  \end{align*}
  belongs to $F(T)_{r-1}$. By induction hypothesis, $S\in F(T)_N + F(J)_{r-1}$.
  Then $t = A(v,w_1,\ldots, w_n) - s \in F(J)_{r} + F(T)_N + F(J)_{r-1} =
    F(T)_N + F(J)_r$.
  This concludes the proof.
\epfv
\begin{cor}
  We have $T = F(T)_N + J$ and $T/J$ is a finitely generated $R$-module.
\end{cor}

\proof
  The first assertion follows immediately from Proposition \ref{pp}; the second
  follows since $F(T)_N$ is finitely generated.
\epfv

In the following we will fix a finite-dimensional associative algebra $P$ and a symmetric linear
function $\phi$ on $P$. Recall that the cenrtal charge of $V$ is $c$. 
For $V$-$P$-bimodules $\tilde W_0\ldots \tilde W_n$ with $\tilde W_0 = \tilde W_n$
and $P$-intertwining operators $\Y_i$ of type $\binom{\tilde W_{i-1}}{W_i \tilde W_i}$,
$i=1\ldots n$, and for $w_i\in W_i$, $i=1,\ldots,n$, we consider the map
\begin{align*}
  &F_{\Y_1,\ldots,\Y_n}(w_1,\ldots,w_n; z_1, \ldots, z_n; q)\\
  &\qquad =\Tr \gYq{1}\cdots \gYq{n}q^{L(0) - \frac{c}{24}}
\end{align*}
and we extend it to an $R$-module map
$\psi_{\Y_1,\ldots,\Y_n} : T \rightarrow G_q((q))[\log q]$ defined by
\begin{align*}
  f\otimes w_1\otimes \ldots \otimes w_n \mapsto
  f\cdot F(w_1,\ldots, w_1; z_1,\ldots,z_n;q).
\end{align*}
We will also use the notation
$F^\phi_{\Y_1,\ldots,\Y_n}$
and
$\psi^\phi_{\Y_1,\ldots,\Y_n}$
when we need to specify the dependence on the symmetric function $\phi$.

\begin{prop}
  The submodule $J$ is contained in the kernel of $\psi_{\Y_1,\ldots,\Y_n}$;
  hence $\psi_{\Y_1,\ldots,\Y_n}$
  induces a map, also denoted by $\psi_{\Y_1,\ldots,\Y_n}$, from $T/J$
  to $G_q((q))[\log q]$.
\end{prop}
\proof
  This follows from applying $\psi_{\Y_1,\ldots, \Y_n}$ to (\ref{rf3}) with
  $l=2$; the resulting equation implies
  $\psi_{\Y_1,\ldots, \Y_n}(A_j(v; w_1,\ldots,w_n))=0$ for all $v\in V$,
  $w_1\in W_1, \ldots, w_n\in W_n$.
\epfv

\begin{prop} For any $w_1\in W_1, \ldots w_n \in W_n$, we have the
  $L(-1)$-derivative property
  \begin{align} \label{FL-1-derivative}
    &\frac{\partial}{\partial z_i} F_{\Y_1,\ldots, \Y_n}
      (w_1,\ldots,w_n; z_1, \ldots, z_n; q)\nn
    &\qquad = F_{\Y_1,\ldots, \Y_n}
      (w_1,\ldots, w_{j-1}, L(-1)w_j, w_{j+1},\ldots,w_n; z_1, \ldots, z_n; q)
  \end{align}
\end{prop}

\proof
  This is an immediate consequence of Lemma \ref{L-1-derivative}.
\epfv

\begin{prop} 
  Consider grading-restricted generalized $V$-modules $W_i$, $\tilde W_i$
  for $i=1,\ldots,n$, with $\tilde W_0 = \tilde W_n$, and logarithmic intertwining
  operators $\Y_i$ of type $\binom{\tilde W_{i-1}}{W_i \tilde W_i}$ for
  $i=1,\ldots,n$. Moreover, suppose $\tilde W_0$ is a $V$-$P$-bimodule
  projective as a right $P$-module.
  Then for any homogeneous elements $w_1\in W_1,\ldots, w_n\in W_n$,
  and any $j=1,\ldots,n$, we have
  \begin{align} \label{O-operator}
    &\left((2\pi i)^2 q \frac{\partial}{\partial q} +
      \tilde G_2(q) \sum_{i=1}^n \wt w_i +
      \tilde G_2(q) \sum_{i=1}^n z_i \frac{\partial}{\partial z_i} -
      \sum_{i\neq j} \tilde \wp_1(z_i - z_j; q)\frac{\partial}{\partial z_i}\right)\cdot\bn
      &\qquad\qquad \cdot F_{\Y_1,\ldots,\Y_n}(w_1,\ldots,w_n; z_1, \ldots ,z_n; q)\bn
      &\qquad\quad + \tilde G_2(q)\sum_{i = 1}^n 
      F_{\Y_2,\ldots,\Y_n}
        (w_1,\ldots,w_{i-1},L(0)_n w_i, w_{i+1},\ldots ,w_n; z_1, \ldots ,z_n; q)\bn
      &\quad= F_{\Y_2,\ldots,\Y_n}(w_1,\ldots, w_{j-1},
        L(-2)w_j, w_{j+1},\ldots,w_n; z_1, \ldots ,z_n; q)\bn
      &\quad\quad - \sum_{k=1}^\infty \tilde G_{2k+2}(q)\cdot \bn
      &\qquad\qquad \cdot F_{\Y_1,\ldots,\Y_n}(w_1,\ldots, w_{j-1},
        L(2k)w_j, w_{j+1},\ldots,w_n; z_1, \ldots ,z_n; q)\bn
      &\quad\quad + \sum_{i\neq j}\sum_{m=1}^\infty (-1)^m 
        \tilde \wp_{m+1}(z_i - z_j); q)\cdot \bn
      &\qquad\qquad \cdot F_{\Y_1,\ldots,\Y_n}(w_1,\ldots, w_{i-1},
        L(m-1)w_i, w_{i+1},\ldots,w_n; z_1, \ldots ,z_n; q)
  \end{align}
\end{prop}
\proof
  This follows from (\ref{rf4}) and the definition of $\psi_{\Y_1,\ldots,\Y_n}$,
  using Lemma \ref{x-L-0-derivative} and (\ref{FL-1-derivative}).
\epfv

\begin{rema} {\rm
  This formula differs from the one obtained in \cite{H2}, by the presence of
  the additional term $$G_2(q)\sum_{i = 1}^n 
      F_{\Y_2,\ldots,\Y_n}
        (w_1,\ldots,w_{i-1},L(0)_n w_i, w_{i+1},\ldots ,w_n; z_1, \ldots ,z_n; q)$$
  due to the non-semisimplicity of the operator $L(0)$.
  }
\end{rema}

We will now consider a second grading on the space $T$ by ``modular weights'':
define the modular weight on the ring $R$ by assigning weight $2k$ to the
function $G_{2k}(q)$ and weight $m$ to $\tilde \wp_m(z; q)$, and denote
by $R_p$ the homogeneous subspace of modular weight $p$.
Now if $t = f\otimes w_1\otimes\ldots\otimes w_n\in T$, with homogeneous
$f\in R_p$ and $w_i\in W_i$, $i=1,\ldots, n$, we assign $t$
modular weight $p + \sum_{i=1}^n \wt w_i$.
Clearly from this definition, if $v\in V$ and $w_i\in W_i$,
$i=1,\ldots,n$ are homogeneous, then the element 
$A(v;w_1,\ldots,w_n) \in J$ has modular weight
$\wt v + \sum_{i=1}^n \wt w_i + 1$. As a consequence, we have

\begin{prop}
  The ideal $J$ is graded by modular weights; in particular, this grading
  induces a grading on the quotient module $T/J$.
\end{prop}

\begin{prop}
  Let $W_1,\ldots W_n$ be grading-restricted generalized modules for the vertex operator algebra
  $V$, and consider homogeneous $w_i\in W_i$ for $i = 1,\ldots n$.
  Then there exist elements $a_{p,i} \in R_p$ for $p=1,\ldots, m$
  such that for any $V$-$P$-bimodules $\tilde W_j$ and $P$-intertwining operators
  $\Y_j$ of type $\binom{\tilde W_{j-1}}{W_j\tilde W_j}$, $j=1,\ldots,n$,
  with $\tilde W_0 = \tilde W_n$ projective as a right $P$-module,
  the series
  $$
    \Tr \gYq{1}\cdot \ldots \cdot \gYq{n} q^{L(0) - \frac{c}{24}}
  $$
  satisfies the differential equations
  \begin{equation}
    \der{m}{z_i} \varphi + \sum_{p=1}^m a_{p,i}(z_1,\ldots,z_n; q) \der{m-p}{z_i} \varphi = 0.
  \end{equation}
  for $i=1,\ldots, n$, in the region $1 > |q_{z_1}| >  \ldots > |q_{z_1}| > |q| > 0$.
\end{prop}

\proof
  Consider the submodule $M_i$ of $T/J$ generated by the elements
  $$
    1\otimes w_1\otimes \ldots \otimes w_{i-1}\otimes L(-1)^kw_i
      \otimes w_{i+1}\otimes \ldots w_n + J
  $$
  for $k\in \N$. Since $R$ is Noetherian and $T/J$ is finitely generated,
  $M_i$ is also finitely generated. Therefore there exists an integer $m$
  and elements $a_{p,i}\in R$ such that
  \begin{align}
    &1\otimes w_1\otimes \ldots \otimes w_{i-1}\otimes L(-1)^mw_i
      \otimes w_{i+1}\otimes \ldots w_n + J\nn
    &\qquad = \sum_{i=1}^m a_{p,i}(z_1,\ldots,z_n; q)\cdot\nn
    &\qquad\qquad \cdot 1\otimes w_1\otimes \ldots \otimes w_{i-1}\otimes L(-1)^{m-p}w_i
      \otimes w_{i+1}\otimes \ldots w_n + J.
  \end{align}
  Note that since the modular weight of 
  $
    1\otimes w_1\otimes \ldots \otimes w_{i-1}\otimes L(-1)^kw_i
    \otimes w_{i+1}\otimes \ldots w_n + J
  $
  is $\sum_{i=1}^n \wt w_i + k$, we can choose the element $a_{p,i}$ to have degree $p$.
  The conclusion follows by applying the map $\psi_{\Y_1,\ldots, \Y_n}$ to both sides
  and using the $L(-1)$ derivative property.
\epfv

Given a logarithmic intertwining operator $\Y$ of type $\binom{W_3}{W_1 W_2}$,
for grading-restricted generalized modules $W_1$, $W_2$, $W_3$, we will consider the map
\begin{align*}
  \Y^{(k)} :& W_1\otimes W_2 \rightarrow W_3\{x\}[\log x]\\
    & w_1\otimes w_2 \mapsto \Y(L(0)_n^kw_1, x)w_2
\end{align*}
for $k\in \N$. Since $L(0)_n$ acts on $W_1$ as a $V$-module endomorphism,
$\Y^{(k)}$ is itself an intertwining operator of the same type.
Also, for $j = 1,\ldots,n$ define 
\begin{align*}
  &Q_j(1\otimes w_1\otimes \ldots \otimes w_n) =\\
      &\qquad w_1\otimes \ldots\otimes w_{j-1}\otimes
        L(-2)w_j\otimes w_{j+1}\otimes\ldots\otimes w_n\bn
      &\quad\quad - \sum_{k=1}^\infty \tilde G_{2k+2}(q)\cdot \\
      &\qquad\qquad \cdot w_1\otimes\ldots\otimes w_{j-1}\otimes
        L(2k)w_j\otimes w_{j+1}\otimes\ldots\otimes w_n\bn
      &\quad\quad + \sum_{i\neq j}\sum_{m=1}^\infty (-1)^m 
        \tilde \wp_{m+1}(z_i - z_j; q)\cdot \bn
      &\qquad\qquad \cdot w_1\otimes\ldots\otimes w_{i-1}\otimes
        L(m-1)w_i\otimes w_{i+1}\otimes \ldots\otimes w_n.
\end{align*}
Then, according to (\ref{O-operator}), we have
\begin{align} \label{rf7}
  &\psi_{\Y_1, \ldots, \Y_n}(Q_j(1\otimes w_1 \otimes \ldots \otimes w_n)) =\bn
  &\qquad \left((2\pi i)^2 q \frac{\partial}{\partial q} +
    \tilde G_2(q) \sum_{i=1}^n \wt w_i +
    \tilde G_2(q) \sum_{i=1}^n z_i \frac{\partial}{\partial z_i} -
    \sum_{i\neq j} \tilde \wp_1(z_i - z_j; q)\frac{\partial}{\partial z_i}\right)\cdot\bn
  &\qquad\qquad \cdot F_{\Y_1,\ldots,\Y_n}(w_1,\ldots, w_n; z_1,\ldots z_n; q)\bn
  &\qquad\quad + \tilde G_2(q)\sum_{i = 1}^n f(q)
    F_{\Y_1,\ldots, \Y_{i-1},\Y_i^{(1)}, \Y_{i+1}, \ldots\Y_n}
      (w_1,\ldots, w_n; z_1, \ldots, z_n; q).
\end{align}
For $\alpha \in \C$ and $j=1,\ldots,n$, we define the differential operator
\begin{align}\label{O-definition}
  &\mathcal O_j(\alpha) = \left((2\pi i)^2 q \frac{\partial}{\partial q} +
    \tilde G_2(q) \alpha +
    \tilde G_2(q) \sum_{i=1}^n z_i \frac{\partial}{\partial z_i} -
    \sum_{i\neq j} \tilde \wp_1(z_i - z_j; q)\frac{\partial}{\partial z_i}\right)
\end{align}
and introduce the notation
\begin{align*}
  &\D_j(\alpha) F_{\Y_1,\ldots, \Y_n}(w_1,\ldots,w_n; z_1,\ldots, z_n; q) \bn
  &\qquad = \sO_j(\alpha)F_{\Y_1,\ldots, \Y_n}(w_1,\ldots,w_n; z_1,\ldots, z_n; q) \bn
  &\qquad\qquad + \tilde G_2(q) \sum_{i=1}^n 
    F_{\Y_1,\ldots, \Y_n}
      (w_1,\ldots, w_{i-1}, L(0)_nw_i, w_{i+1},\ldots ,w_n; z_1,\ldots, z_n; q)
\end{align*}
and inductively
\begin{align*}
  &\left(\prod_{l=1}^k \D_j(\alpha_l)\right)
    F_{\Y_1,\ldots, \Y_n}(w_1,\ldots,w_n; z_1,\ldots, z_n; q)\bn
    &\qquad = \sO_j(\alpha_1)\left( \left(\prod_{l=2}^k \D_j(\alpha_l)\right)
    F_{\Y_1,\ldots, \Y_n}(w_1,\ldots,w_n; z_1,\ldots, z_n; q)\right)\bn
  &\qquad \qquad + \tilde G_2(q)\sum_{i=1}^n
    \left(\prod_{l=2}^k \D_j(\alpha_l)\right)\cdot\bn
  &\qquad\qquad \qquad \cdot F_{\Y_1,\ldots, \Y_n}
      (w_1,\ldots, w_{i-1}, L(0)_nw_i, w_{i+1},\ldots ,w_n; z_1,\ldots, z_n; q).
\end{align*}
\begin{rema}{\rm
  Using this notation, (\ref{rf7}) can be written as
  \begin{equation*}
    \psi_{\Y_1, \ldots, \Y_n}(\sQ_j(1\otimes w_1\otimes \ldots \otimes w_n)) = 
      \D_j\left(\sum_{i=1}^n \wt w_i\right) \psi_
      {\Y_1, \ldots, \Y_n}(1\otimes w_1\otimes \ldots \otimes w_n).
  \end{equation*}
  We now want to extend $Q_j$ to a map $\sQ_j: T\rightarrow T$ such that, for any
  $t\in T$ of modular weight $\alpha$,
  \begin{align*}
    &\psi_{\Y_1, \ldots, \Y_n}(\sQ_j(t)) = \D_j(\alpha)\psi_
      {\Y_1, \ldots, \Y_n}(t).
  \end{align*}
  This is necessary since we will need to apply $\psi_{\Y_1,\ldots, \Y_n}$ to repeated
  iterations of the map $\sQ_j$ on the element $1\otimes w_1\otimes \ldots
  \otimes w_n$ for homogeneous elements $w_1,\ldots, w_n$;
  and while the elements of $T$ obtained this way are not homogeneous in
  the conformal grading, they have a well defined modular weight.
  }
\end{rema}

\begin{defn}{\rm
    We define functions $\vartheta_j: R\rightarrow R$ for $j=1,\ldots, n$
    in the following way: let
    $$
    \vartheta_j(\tilde G_k(q)) =
      (2\pi i)^2 q\der{}{q}\tilde G_k(q) + k\tilde G_2(q)\tilde G_k(q),
    $$
    and for the formal series  $\tilde \wp_m(z_r - z_s; q_\tau)$ with
    $1\leq r \neq s\leq n$ and any $m\geq 2$, and $j = 1,\ldots, n$, we
    define $\vartheta_j(\tilde \wp_m)$ by
    \begin{align*}
      &\vartheta_j(\tilde\wp_m(z_r - z_s; q))\bn
      &\qquad = (2\pi i)^2 q \der{}{q}\tilde \wp_m(z_r-z_s;q) 
        + m\tilde G_2(q)\tilde\wp_m(z_r-z_s;q)\bn
      &\qquad \quad - m\tilde G_2(q)(z_r - z_s)
        \tilde\wp_{m+1}(z_r-z_s;q)\bn
      &\qquad \quad + m\tilde\wp_{m+1}(z_r-z_s; q)
        \left(\tilde\wp_1(z_r-z_j;q) - \tilde\wp_1(z_s-z_j;q)\right)
    \end{align*}
    if $j\notin \{r,s\}$; and by
    \begin{align*}
      &\vartheta_j(\tilde \wp_m(z_j - z_s; q))\bn
      &\qquad = (2\pi i)^2 q \der{}{q}\tilde \wp_m(z_j-z_s;q)
        + m\tilde G_2(q)\tilde \wp_m(z_j-z_s;q)\bn
      &\qquad \quad - m\tilde G_2(q)(z_j - z_s)\tilde \wp_{m+1}(z_j-z_s;q)\bn
      &\qquad \quad - m\tilde \wp_{m+1}(z_j-z_s; q)
        \tilde \wp_1(z_s-z_j;q)
    \end{align*}
    if $j = r$, and we extend $\vartheta_j$ as a derivation on the ring $R$.
  }
\end{defn}

\begin{prop}
  Let $\varphi(z_1,\ldots,z_n; q) = 
  \vartheta_j(\tilde\wp_m(z_r - z_s; q))$; then $\varphi(z_1,\ldots, z_n; q_\tau)$
  converges uniformly to an elliptic function in the variables $z_1,\ldots, z_n$
  with possible poles at $z_r = z_s + m\tau + n$, $n, m\in \Z$, $r,s=1,\ldots,n$.
  Moreover, for any $g\in SL_2(\Z)$, if
  $$
    g = \left(
      \begin{array}{c c}
        \alpha & \beta\\
        \gamma & \delta
      \end{array}
    \right)
  $$
  we have
  \begin{align*}
    \varphi\left(\frac{z_1}{\gamma\tau + \delta},\ldots,\frac{z_n}{\gamma\tau + \delta};
    \frac{\alpha \tau + \beta}{\gamma\tau + \delta}\right) = 
    \left(
      \gamma\tau + \delta
    \right)^{m+2} \varphi(z_1,\ldots, z_n; \tau).
  \end{align*}
\end{prop}
\proof
  Easy computation using transformation properties of $\wp_m$ and $G_k$ under
  the action of $SL_2(\Z)$.
\epfv

\noindent 
In particular, $\vartheta_j$ is indeed a map from $R$ to $R$ and
if $f\in R$ has modular weight $p$, then $\vartheta_j(f)$ has modular weight $p+2$.
We then consider the $\C$-linear maps $\sQ_j : T\rightarrow T$ for $j=1,\ldots, n$
defined by
\begin{align*}
  &\sQ_j(f(q)\otimes w_1\otimes \ldots \otimes w_n) \bn
  &\qquad = f(q)\cdot Q_j(1\otimes w_1\otimes \ldots \otimes w_n)
     + \vartheta_j(f(q))\otimes w_1\otimes \ldots \otimes w_n;
\end{align*}
note that if the modular weight of $t\in T$ is $\alpha$, then
the modular weight of $\sQ_j(t)$ is $\alpha + 2$.

\begin{prop}
  Let $t\in T$ be an element of modular weight $\alpha$. Then
  \begin{align*}
    \psi_{\Y_1,\ldots,\Y_n}(\sQ_j(t)) = \D_j(\alpha)\cdot\psi_{\Y_1,\ldots, \Y_n}(t).
  \end{align*}
\end{prop}
\proof
  It is enough to prove this for elements $t$ of the form
  $\tilde G_k(q)\otimes w_1\otimes \ldots\otimes w_n$ or
  $\tilde \wp_m(z_i - z_j; q)\otimes w_1\otimes \ldots\otimes w_n$.
  Suppose $t$ is of the first form with $\sum_{i=1}^n \wt w_i = s$
  with $k + s = \alpha$.
  Then
  \begin{align*}
    &\psi_{\Y_1,\ldots,\Y_n}(\sQ_j(t)) \nn
    &\qquad = \psi_{\Y_1,\ldots,\Y_n}(\tilde G_k(q) 
      Q_j(1\otimes w_1\otimes\ldots\otimes w_n) +
      \theta(\tilde G_k(q))\otimes w_1\otimes\ldots\otimes w_n)\bn
    &\qquad = \tilde G_k(q)\D_j(s)F_{\Y_1,\ldots,\Y_n}
      (w_1,\ldots, w_n; z_1, \ldots, z_n; q)\bn
    &\qquad \quad + \left((2\pi i)^2 q \der{}{q}\tilde G_k(q) + 
      k \tilde G_2(q)\tilde G_k(q)\right)F_{\Y_1,\ldots,\Y_n}
      (w_1,\ldots, w_n; z_1, \ldots, z_n; q)\bn
    &\qquad = \D_j(s+k)\tilde G_k(q)F_{\Y_1,\ldots,\Y_n}
      (w_1,\ldots, w_n; z_1, \ldots, z_n; q)\bn
    &\qquad = \D_j(\alpha)\psi_{\Y_1,\ldots, \Y_n}(t).
  \end{align*}
  The proof of the other case is a similar computation.
\epfv
\begin{prop} 
  Let $\alpha = \sum_{i=1}^n \wt w_i$; then for any $s\in \N$, we have
  \begin{align}\label{rf5}
    &\psi_{\Y_1,\ldots,\Y_n}(Q_j^s(1\otimes w_1 \otimes \ldots \otimes w_n))\bn
    &\qquad = \prod_{i=1}^s \D_j(\alpha + 2(s - i)) \cdot 
      F_{\Y_1,\ldots, \Y_n}(w_1,\ldots,w_n; z_1,\ldots, z_n; q).
  \end{align}
\end{prop}
\proof
  We proceed by induction on $s$; the base case $s = 0$ follows by definition of
  $\psi_{\Y_1,\ldots, \Y_n}$. Now suppose the claim holds for $s-1$; then
  by (\ref{O-operator}), and since $Q^k(1\otimes w_1\otimes\ldots\otimes w_n)$
  has modular weight equal to
  $\alpha + 2(k-1)$,
  \begin{align*}
    &\psi_{\Y_1,\ldots,\Y_n}(Q_j^s(1\otimes w_1 \otimes \ldots \otimes w_n))\bn
    &\qquad = \sO_j(\alpha + 2(s-1)) \psi_{\Y_1,\ldots,\Y_n}
      (Q_j^{s-1}(1\otimes w_1 \otimes \ldots \otimes w_n))\bn
    &\qquad\qquad + \tilde G_2(q) \sum_{i=1}^n
      \psi_{\Y_1,\ldots,\Y_{i-1},\Y_i^{(1)},\Y_{i+1},\ldots,\Y_n}
        (Q_j^{s - 1}(1\otimes w_1 \otimes \ldots \otimes w_n))\bn
    &\qquad = \sO_j(\alpha + 2(s-1)) \left(\prod_{l=1}^{s-1}
      \D_j(\alpha + 2(s - 1 - l))\right)\cdot \bn
    &\qquad\qquad\qquad
      \cdot F_{\Y_1,\ldots, \Y_n}(w_1,\ldots,w_n; z_1,\ldots, z_n; q)\bn
    &\qquad\qquad + \tilde G_2(q) \sum_{i=1}^n
      \left(\prod_{l=1}^{s-1}
      \D_j(\alpha + 2(s - 1 - l))\right)\cdot \bn
    &\qquad\qquad\qquad
      \cdot F_{\Y_1,\ldots, \Y_n}
      (w_1,\ldots,w_{i-1}, L(0)_nw_i, w_{i+1},\ldots ,w_n; z_1,\ldots, z_n; q)\bn
    &\qquad = \sO_j(\alpha + 2(s-1)) \left(\prod_{l=2}^{s}
      \D_j(\alpha + 2(s - l))\right)\cdot \bn
    &\qquad\qquad\qquad
      \cdot F_{\Y_1,\ldots, \Y_n}(w_1,\ldots,w_n; z_1,\ldots, z_n; q)\bn
    &\qquad\qquad + \tilde G_2(q) \sum_{i=1}^n
      \left(\prod_{l=2}^{s}
      \D_j(\alpha + 2(s - l))\right)\cdot \bn
    &\qquad\qquad\qquad
      \cdot F_{\Y_1,\ldots, \Y_n}
      (w_1,\ldots,w_{i-1}, L(0)_nw_i, w_{i+1},\ldots ,w_n; z_1,\ldots, z_n; q)\bn
    &\qquad = \prod_{l=1}^s \D_j(\alpha + 2(s - l)) \cdot 
      F_{\Y_1,\ldots, \Y_n}(w_1,\ldots,w_n; z_1,\ldots, z_n; q),
  \end{align*}
  which concludes the proof.
\epfv

\begin{prop} \label{diffeqprop}
  Let $W_1,\ldots W_n$ be grading-restricted generalized modules for the vertex operator algebra
  $V$, and consider homogeneous $w_i\in W_i$ for $i = 1,\ldots n$; let
  $\alpha = \sum_{i=1}^n \wt w_i$.
  Then there exist elements $b_{p,i} \in R_{2p}$ for $p=1,\ldots, m$
  such that for any $V$-$P$-bimodules $\tilde W_j$ and $P$-intertwining operators
  $\Y_j$ of type $\binom{\tilde W_{j-1}}{W_j\tilde W_j}$, $j=1,\ldots,n$,
  with $\tilde W_0 = \tilde W_n$ projective as a right $P$-module,
  the series
  $$
    \varphi = F_{\Y_1,\ldots, \Y_n}(w_1,\ldots, w_n; z_1,\ldots z_n; q)
  $$
  satisfies the differential equations
  \begin{align} \label{diffeq}
    &\prod_{l=1}^m \D_j(\alpha + 2(m - l)) 
      \varphi\bn 
    &\qquad + \sum_{p=1}^{m} b_{p,j}(z_1,\ldots,z_n; q)
      \prod_{l=1}^{m-p} \D_j(\alpha + 2(m-p-l))
      \varphi = 0
  \end{align}
  for $j=1,\ldots, n$, in the region $1 > |q_{z_1}| >  \ldots > |q_{z_1}| > |q| > 0$.\\
  Moreover, for any $k=1,\ldots, n$, the series
  \begin{align*}
    &\varphi_k(z_1,\ldots,z_n; q) \bn
    &\qquad = F^\phi_{\Y_1,\ldots, \Y_n}
    (w_1,\ldots, w_{k-1}, L(0)_n w_k, w_{k+1}, \ldots, w_n; z_1,\ldots z_n; q)
  \end{align*}
  also satisfies (\ref{diffeq}), for the same choice of elements $b_{p,i}$.
\end{prop}
\proof
Consider the $R$-submodule $M$ of $T/J$ generated by the elements 
$$
  \sQ_j^k(1\otimes w_1\otimes \ldots \otimes w_n) + J
$$ 
for $k\in \N$. Since $T/J$ is finitely generated and $R$ is noetherian,
$M$ is also finitely generated; therefore, there exists $m$ and
elements $b_{p,j}(z_1,\ldots,z_n; q)\in R$, $p=1,\ldots, m$
such that
\begin{align*}
  &\sQ_j^m(1\otimes w_1\otimes \ldots \otimes w_n)\bn
  &\qquad + \sum_{p=1}^{m} b_{p,j}(z_1,\ldots, z_n; q)
  \sQ_j^{m-p}(1\otimes w_1\otimes \ldots \otimes w_n) \in J.
\end{align*}
Since the modular weight of $\sQ_j^i(1\otimes w_1\otimes \ldots \otimes w_n)$
is $\alpha + 2i$, we can choose $b_{p,j}$ in $R_{2p}$; then applying
$\psi_{\Y_1,\ldots,\Y_n}$ to both sides of the equation and applying
(\ref{rf5}), we obtain (\ref{diffeq}).\\
The last part of the proposition follows by applying the first part to
the intertwining operators $\Y_1,\ldots, \Y_{k-1}, \Y_k^{(1)}, \Y_{k+1},
  \ldots, Y_n$.
\epfv

\begin{rema} {\rm
    Note that the coefficients $b_{p,j}$ only depend on the elements
    $w_1,\ldots, w_n$; in particular, they do not depend on the
    choice of $P$ and $\phi$ in the
    definition of the pseudotrace.
  }
\end{rema}
\begin{rema} \label{rm1}
  {\rm
    For $i_j \in \N$, and $i=1,\ldots, n$, let
    \begin{align*}
      &\varphi_{i_1,\ldots, i_n}(z_1,\ldots, z_n; q)\bn
      &\qquad = F_{\Y_1,\ldots, \Y_n} (L(0)_n^{i_1}w_1,
    \ldots, L(0)_n^{i_n}w_n; z_1,\ldots,z_n; q), 
    \end{align*}
    The differential equation (\ref{diffeq}) depends on all the functions
    $\varphi_{i_1,\ldots, i_n}$; therefore, we obtain a system of differential
    equations of which $\{\varphi_{i_1,\ldots, i_n} | i_j \in \N, j = 1,\ldots,n\}$
    is a solution. Since $L(0)_n$ is locally nilpotent, this system of equations
    is finite; if the modules $W_1,\ldots, W_n$ all have length smaller than $l$,
    then we obtain a system for the $l(n+1)$ functions 
    $\{\varphi_{i_1,\ldots, i_n} | i_j = 0,\ldots, l, j = 1,\ldots,n\}$.
    If the modules are ordinary modules the system decouples and we obtain the
    equations in \cite{H2}.
  }
\end{rema}

\begin{prop}
  The series
  \begin{align} \label{g1s}
    F_{\Y_1,\ldots, \Y_n}(w_1,\ldots, w_n; z_1,\ldots, z_n; q_\tau)
  \end{align}
  is absolutely convergent in the region
  $1>|q_{z_1}| > \ldots > |q_{z_n}| > |q_\tau| > 0|$ and can be extended
  to a multivalued analytic function in the region $\Im (\tau) >0$,
  $z_i \neq z_j + l + m\tau$ for $i\neq j$, $l,m\in \Z$.
\end{prop}
\proof
  For fixed $z_1,\ldots, z_n$ such that $|q_{z_1}| > \ldots >|q_{z_n}| > 0$,
  the coefficients in the variable $q$, $\log q$ of the series
  $$
    F_{\Y_1,\ldots, \Y_n}(w_1,\ldots, w_n; z_1,\ldots, z_n; q)
  $$
  are absolutely convergent, and the series satisfies a system of
  differential equations with a regular singular point at $q=0$ and
  analytic coefficients.
  Since the coefficients of the differential equations are analytic
  functions in $z_1,\ldots, z_n, q_\tau$, with possible singularities
  in the region $\Im (\tau) > 0$, $z_i\neq z_j + l + m\tau$ for
  $i\neq j$ and $l, m \in \Z$, the solutions of the system (\ref{diffeq})
  can be extended to an analytic (multivalued) function in the same region.
\epfv

We will call \emph{genus-one correlation functions} the analytic extensions
of (\ref{g1s}) to the region $\Im (\tau) >0$,
$z_i \neq z_j + l + m\tau$ for $i\neq j$, $l,m\in \Z$,
and we will denote them by
$$
  \overline F^\phi_{\Y_1,\ldots, \Y_n}(w_1,\ldots, w_n; z_1,\ldots, z_n; q_\tau).
$$

\begin{prop}[Genus-one commutativity]
  Let $W_i$ be grading-restricted generalized $V$-modules, $\tilde W_i$ be $V$-$P$ bimodules
  and $\Y_i$ logarithmic intertwining operators of type
  $\binom{\tilde W_{i-1}}{W_i \tilde W_{i-1}}$ for $i=1,\ldots, n$,
  with $\tilde W_0 = \tilde W_n$ projective as a right $P$-module.
  Then for any $k\leq n-1$, there exists a $V$-$P$-bimodules $\hat W_k$,
  and logarithmic intertwining operators $\hat \Y_k$, $\hat \Y_{k+1}$
  of type $\binom{\hat W_k}{W_k\tilde W_{k+1}}$,
  $\binom{\tilde W_{k-1}}{W_{k+1} \hat W_k}$ such that
  \begin{align*}
    &\overline
      F^\phi_{\Y_1, \ldots, \Y_n}(w_1,\ldots, w_n; z_1, \ldots, z_n; \tau)\\
    &\qquad  = \overline 
      F^\phi_{\Y_1, \ldots, \Y_{k-1},\hat \Y_{k+1},\hat Y_k, \Y_{k+2}
      ,\ldots\Y_n}(w_1,\ldots, w_{k-1}, w_{k+1}, w_k, w_{k+2} \ldots, w_n;\\
      &\hspace{7cm}
      z_1,\ldots, z_{i-1}, z_{i+1}, z_i, z_{i+2}, \ldots, z_n; \tau)
  \end{align*}
  as multivalued analytic functions.
\end{prop}
\proof 
  Follows from commutativity for $P$ intertwining operators.
\epfv

\begin{prop}[Genus-one associativity]
  Let $W_i$ be grading-restricted generalized $V$-modules, $\tilde W_i$ be $V$-$P$-bimodules
  and $\Y_i$ logarithmic intertwining operators of type
  $\binom{\tilde W_{i-1}}{W_i \tilde W_{i-1}}$ for $i=1,\ldots, n$,
  with $\tilde W_0 = \tilde W_n$ projective as a right $P$-module.
  Then for any $k\leq n-1$, there exists a $V$-$P$-bimodules $\hat W_k$,
  and logarithmic intertwining operators $\hat \Y_k$, $\hat \Y_{k+1}$
  of type $\binom{\hat W_k}{W_k W_{k+1}}$,
  $\binom{\tilde W_{k-1}}{\hat W_k \tilde W_{k+1}}$ such that the series
  \begin{align*}
    &\overline
    F^\phi_{\Y_1, \ldots, \Y_{k-1}, \hat\Y_{k+1},\Y_{k+2},\ldots, \Y_n}
      (w_1,\ldots, w_{k-1}, \hat \Y(w_k,z_k - z_{k+1})w_{k+1}, \\
      &\hspace{6cm} w_{k+2}, \ldots, w_n;z_1, \ldots, z_n; \tau)\\
    &\qquad  = \sum_{r\in \R}\overline F^\phi_{\Y_1, \ldots, \Y_{k-1},
      \hat\Y_{k+1},\Y_{k+2},\ldots, \Y_n}
      (w_1,\ldots, w_{k-1}, P_r(\hat \Y(w_k,z_k - z_{k+1})w_{k+1}), \\
      &\hspace{6cm} w_{k+2}, \ldots, w_n;z_1, \ldots, z_n; \tau)
  \end{align*}
  is absolutely convergent in the region
  $$
    1 > |q_{z_1}| > \ldots |q_{z_{k-1}}| > |q_{z_{k+1}}|
    \ldots > |q_{z_{n}}| > |q_{\tau}| > 0
  $$
  and $ 1 > |q_{(z_{k} - z_{k+1})}| > 0$ and converges to
  $\overline F^\phi_{\Y_1, \ldots, \Y_n}(w_1,\ldots, w_n;
  z_1, \ldots, z_n; \tau)$ when $1 > |q_{z_1}| > 
    \ldots > |q_{z_{n}}| > |q_{\tau}| > 0$ and 
    $|q_{(z_{k} - z_{k+1})}| > 1 > |q_{(z_{k} - z_{k+1})}-1| > 0$.
\end{prop}
\proof 
  By associativity for $P$-intertwining operators, there exist a
  $V$-$P$-bimodule $\tilde W_k$ and  logarithmic intertwining operators 
  $\hat \Y_k$, $\hat \Y_{k+1}$
  of type $\binom{\hat W_k}{W_k W_{k+1}}$,
  $\binom{\tilde W_{k-1}}{\hat W_k \tilde W_{k+1}}$ such that
  for any $z_1,\ldots, z_n \in \C$ satisfying $1> |q_{z_1}| >
  \ldots > |q_{z_n}| > 0$ and
  $|q_{z_{k+1}}| > |q_{z_k} - q_{z_{k+1}}| > 0$,
  and for any element $\tilde w_n^\prime\in \tilde W_n^\prime$,
  $\tilde w_n\in \tilde W_n$
  \begin{align*}
    &\langle w_n^\prime, \gYq{1}\ldots \gYq{k-1}\cdot\\
    & \qquad \qquad \cdot
      \hat \Y_{k+1}(\U(q_{z_{k+1}})
      \hat \Y_k(w_k, z_k - z_{k+1})w_{k+1}, q_{z_{k+1}})\cdot\\
    &\qquad\qquad \cdot  \gYq{k+2}\ldots \gYq{n}\tilde w_n\\
    &\qquad = \langle \tilde w_n^\prime, \gYq{1}\ldots \gYq{n}
      \tilde w_n \rangle
  \end{align*}
  and therefore as series in $q$ and $\log q$,
  \begin{align*}
    &\Tr\gYq{1}\ldots \gYq{k-1}\cdot\\
    & \qquad \qquad \cdot
      \hat \Y_{k+1}(\U(q_{z_{k+1}})
      \hat \Y_k(w_k, z_k - z_{k+1})w_{k+1}, q_{z_{k+1}})\cdot\\
    &\qquad\qquad \cdot  \gYq{k+2}\ldots \gYq{n}q^{L(0)-\frac{c}{24}}\\
    &\qquad = \Tr \gYq{1}\ldots \gYq{n}q^{L(0) - \frac{c}{24}}. 
  \end{align*}
  The right hand side is absolutely convergent when $q = q_\tau$ and
  $1>|q_{z_1}| > \ldots > |q_{z_1}| > |q_\tau| > 0$, the right hand side is
  also absolutely convergent if $|q_{z_{k}}| > |q_{z_{k}} - q_{z_{k+1}}| > 0$,
  and satisfies the same system of differential equations as
  $$
    \Tr \gYq{1}\ldots \gYq{n}q^{L(0) - \frac{c}{24}}.
  $$
  So in this region the left hand side converges absolutely to a
  function that can be extended to the multivalued analytic function
  $$
    \overline F^\phi_{\Y_1,\ldots,\Y_n}(w_1,\ldots, w_n; z_1,\ldots, z_n, \tau),
  $$
  which concludes the proof.
\epfv

\subsection{Modular invariance of the space of solutions}
\label{s-diffeq-m-inv}

In this section we consider a space of functions which contains the solutions
of the system of differential equations (\ref{diffeq}), and we define an
action of the group $SL_2(\Z)$ on the elements of this space. We prove
that the space of solutions of (\ref{diffeq}) is invariant under
this action.

We will introduce the following notations: let $\chi$ be the space of sequences
(indexed by $n$ indices)
of analytic multivalued functions in the variables $z_1,\ldots, z_n, \tau$,
on the region $\Im(\tau) > 0$, $z_i\neq z_j + n\tau + m$ for
$n,m\in \N$, $i\neq j$
\begin{equation*}
  \Phi = \Phi(z_1,\ldots, z_n; \tau)
    = (\phi_{i_1,\ldots, i_n}(z_1,\ldots , z_n; \tau))_{i_1,\ldots, i_n \in \N}
\end{equation*}
with preferred branches on the region 
$1 > |q_{z_1}| > \ldots > |q_{z_n}| > |q_\tau| > 0$,
such that $\phi_{i_1,\ldots, i_n} \equiv 0$ whenever $\max\{i_1,\ldots, i_n\}$ is
sufficiently large. We will also denote an element $\Phi$ of $\chi$ by
$(\phi_{i_1,\ldots, i_n})$ or using a multi-index notation $(\phi_\mu)$ for
$\mu$ ranging over $\N^n$.
The sum of sequences of this kind is defined component by component, so
that if $\Phi^1$, $\Phi^2$ are two elements of $\chi$, the $\mu$-th component
of $\Phi^1 + \Phi^2$ is
  $(\phi^1_{\mu} + \phi^2_{\mu})$
and similarly we can define the product by another analytic function. Moreover, we
extend differential operators to $\chi$ component-wise:
\begin{equation*}
  \left(\der{}{z_i} \Phi\right)_\mu = \frac{\partial \phi_{\mu}}
  {\partial z_i},\ i=1,\ldots, n.
\end{equation*}
For $j=1,\ldots,n$, let $d_j: \chi \rightarrow \chi$ be the shift operator on
the $j$-th coordinate defined by
\begin{equation*}
  (d_j \Phi)_{i_1,\ldots, i_n} = \phi_{i_1,\ldots, i_{j-1}, i_j+1, i_{j+1},
    \ldots i_n}.
\end{equation*}
Note that for any $\Phi\in \chi$, $d_j^k \Phi = 0$ if $k$ is large enough; therefore
for any function $f(z_1,\ldots, z_n; \tau)$, the operator
\begin{align*}
  e^{f(z_1,\ldots, z_n; \tau) d_j} = 
    \sum_{k=0}^\infty \frac{f^k(z_1,\ldots, z_n; \tau)}{k!}d_j^k
\end{align*}
is well defined.

Now, for $\alpha \in \C$ and $j= 1,\ldots n$, we define
$\D_j(\alpha) : \chi \rightarrow \chi$ by
\begin{align*}
  \D_j(\alpha) = \sO_j(\alpha) + G_2(\tau)\sum_{i=1}^n d_i
\end{align*}
where $\sO_j(\alpha)$ is defined as in (\ref{O-definition}) 
(with $(2\pi i)^2q\der{}{q} = 2\pi i\der{}{\tau}$) and extended component-wise
to $\chi$.

\begin{rema}{\rm
  Using the notation as in remark (\ref{rm1}), for fixed elements $w_i\in W_i$,
  $i=1,\ldots, n$, and logarithmic intertwining operators $\Y_1,\ldots \Y_n$,
  one can consider the element of $\chi$:
  \begin{equation*}
    \Phi = (\varphi_{i_1,\ldots, i_n})_{i_1,\ldots, i_n \in \N}.
  \end{equation*}
  Then if $b_{p,j}(z_1,\ldots, z_n; \tau)$ $j=1,\ldots,n$ are defined as in
  Proposition \ref{diffeqprop}, by the same proposition we have 
  \begin{align}\label{system}
    &\left(\prod_{l=1}^m \D_j(\alpha + 2(m - l)) \right.
      \nn 
    &\left.\qquad + \sum_{p=1}^{m} b_{p,j}(z_1,\ldots,z_n; q)
      \prod_{l=1}^{m-p} \D_j(\alpha + 2(m-p-l))\right)\Phi
      = 0.
  \end{align}
  }
\end{rema}

\begin{defn}[$SL_2(\Z)$ action]{\rm Let $a \in \C$, and consider an element
  $g$ in $SL_2(\Z)$,
  $$
    g = \left(
      \begin{array}{c c}
        \alpha & \beta\\
        \gamma & \delta
      \end{array}
    \right).
  $$
  For $\Phi \in \chi$, we define
  \begin{align} \label{sl2-action}
    &\Phi|_{g,a}(z_1,\ldots,z_n; \tau)\bn
    &\qquad = \left(\frac{1}{\gamma \tau + \delta}\right)^a
    \prod_{i=1}^n e^{-\log(\gamma \tau + \delta)d_i}
    \Phi\left(
      \frac{z_1}{\gamma \tau  + \delta}, \ldots, \frac{z_n}{\gamma\tau + \delta};
      \frac{\alpha \tau + \beta}{\gamma\tau + \delta}\right)
  \end{align}
}
\end{defn}

\begin{prop}
  For any $a\in \C$, (\ref{sl2-action}) defines an action of the group
  $SL_2(\Z)$ on the space $\chi$.
\end{prop}
We will also use the notation $z^\prime$, $\tau^\prime$ to denote
$\frac{z}{\gamma\tau +\delta}$ and
$\frac{\alpha \tau + \beta}{\gamma \tau + \delta}$ respectively.
Note that for any function $f(z_1,\ldots, z_n; q)$, $g\in SL_2(\Z)$,
and any $\Phi\in \chi$,
$$
  (f\Phi)|_{g,a} = f(z_1^\prime, \ldots, z_n^\prime, \tau^\prime)\Phi|_{g,a};
$$
in particular, if $f\in R_p$, then 
$f(z_1^\prime, \ldots, z_n^\prime; \tau^\prime) = 
\left(\gamma \tau + \delta\right)^pf(z_1,\ldots, z_n; \tau)$ and thus 
\begin{equation*}
  (f\Phi)|_{g,a} = f(z_1, \ldots, z_n, \tau)\Phi|_{g,a-p}.
\end{equation*}

\begin{prop}\label{diffeq-invariance}
  Let $\Phi\in \chi$, $a\in \C$ and $g\in G$. Then
  \begin{equation}\label{D-shift}
     \D_j(a)(\Phi|_{g,a}) = (\D_j(a)\Phi)|_{g,a+2}.
  \end{equation}
\end{prop}
\proof
  This is just a straightforward computation, using the transformation properties
  of the functions $G_2(\tau)$ and $\wp_1(z; \tau)$:
  for simplicity we use the notation
  $$
    e^d = \prod_{i=1}^ne^{-\log(\gamma\tau + \delta)d_i}
  $$
  \begin{align*}
    &\D_j(a)(\Phi|_{g,a}) = \D_j(a)
      \left(\left(\frac{1}{\gamma\tau + \delta}\right)^a
      e^d \Phi(z_1^\prime,\ldots, z_n^\prime;\tau^\prime)\right) \bn
    &\qquad = \left((2\pi i) \frac{\partial}{\partial \tau} +
      G_2(\tau) \left(a + \sum_{i=1}^n d_i\right) +
      G_2(\tau) \sum_{i=1}^n z_i \frac{\partial}{\partial z_i} -
      \sum_{i\neq j} \wp_1(z_i - z_j; \tau)\frac{\partial}{\partial z_i}\right)\cdot \bn
    &\qquad \qquad \cdot \left(\left(\frac{1}{\gamma\tau + \delta}\right)^a
      \prod_{i=1}^ne^{-\log(\gamma \tau + \delta) d_i}
      \Phi(z_1^\prime\ldots, z_n^\prime;\tau^\prime)\right)\bn
    &\qquad = -(2\pi i) \gamma a \left(\frac{1}{\gamma\tau + \delta}\right)^{a+1}
      e^d \Phi(z_1^\prime,\ldots, z_n^\prime, \tau^\prime) \bn
    &\qquad \quad -(2\pi i) \gamma \left(\frac{1}{\gamma\tau + \delta}\right)^{a+1}
      e^d\sum_{i=1}^n d_i \Phi(z_1^\prime,\ldots, z_n^\prime \tau^\prime) \bn
    &\qquad \quad -(2\pi i) \gamma \left(\frac{1}{\gamma\tau + \delta}\right)^{a+1}
      e^d\sum_{i=1}^n z_i^\prime
      \frac{\partial \Phi}{\partial z_i}(z_1^\prime,\ldots, z_n^\prime, \tau^\prime) \bn
    &\qquad \quad +(2\pi i) \left(\frac{1}{\gamma\tau + \delta}\right)^{a+2} e^d      
      \frac{\partial \Phi}{\partial \tau}(z_1^\prime\ldots,z_n^\prime, \tau^\prime) \bn
    &\qquad \quad + \left(G_2(\tau^\prime)\left(\frac{1}{\gamma\tau + \delta}\right)^2
      + 2\pi i \gamma \left(\frac{1}{\gamma\tau + \delta}\right)\right)
      \left(a + \sum_{i=1}^n d_i\right)\cdot\bn
    &\qquad\qquad\qquad \cdot\left(\frac{1}{\gamma\tau + \delta}\right)^a
      e^d\Phi(z_1^\prime,\ldots, z_n^\prime; \tau^\prime) \bn
    &\qquad \quad + \left(G_2(\tau^\prime) \left(\frac{1}{\gamma\tau + \delta}\right)^2
      + 2\pi i \gamma \left(\frac{1}{\gamma\tau + \delta}\right)\right)\cdot\bn
    &\qquad\qquad\qquad\cdot
      \left(\frac{1}{\gamma\tau + \delta}\right)^{a} e^d\sum_{i=1}^n  z_i^\prime
      \frac{\partial \Phi}{\partial z_i}(z_1^\prime,\ldots, z_n^\prime, \tau^\prime)\bn
    &\qquad \quad - \left(\frac{1}{\gamma\tau + \delta}\right)^{a+2}e^d \sum_{i\neq j}
      \wp_1(z_i^\prime - z_j^\prime; \tau^\prime)
      \frac{\partial \Phi}{\partial z_i}(z_1^\prime,\ldots, z_n^\prime, \tau^\prime)\bn
    &\qquad = \left(\frac{1}{\gamma\tau + \delta}\right)^{a+2}e^d\cdot\bn
    &\qquad \qquad \cdot \Bigg((2\pi i)
      \frac{\partial \Phi}{\partial \tau}(z_1^\prime,\ldots, z_n^\prime, \tau^\prime)
     + G_2(\tau^\prime)
      \left(a + \sum_{i=1}^n d_i\right)\Phi(z_1^\prime,\ldots, z_n^\prime; \tau^\prime)\\
    &\qquad\qquad\qquad + G_2(\tau^\prime)
      \sum_{i=1}^n  z_i^\prime
      \frac{\partial \Phi}{\partial z_i}(z_1^\prime,\ldots, z_n^\prime, \tau^\prime)
      - \sum_{i\neq j}
      \wp_1(z_i^\prime - z_j^\prime; \tau^\prime)
      \frac{\partial \Phi}{\partial z_i}(z_1^\prime,\ldots, z_n^\prime, \tau^\prime)\Bigg)
  \end{align*}
  which is equal to $(\D_j(a)\Phi)|_{g,a+2}$, concluding the proof.
\epfv

We can then prove the following

\begin{prop}\label{diffeq-mod-inv}
  Let $\Phi$ be a solution of the system of differential equations (\ref{system}).
  Then for any $g\in SL_2(\Z)$, $\Phi|_{g, \alpha}$ is also a solution
  of the same system.
\end{prop}
\proof
  Just apply $|_{g, \alpha + 2m}$ to both sides of (\ref{system}); since
  $b_{p,j}$ belongs to $R_{2p}$, we find
  \begin{align*}
    &\left(b_{p,j}(z_1,\ldots,z_n; q)
      \prod_{l=1}^{m-p} \D_j(\alpha + 2(m-p-l))\Phi\right)
      \Bigg|_{g, \alpha + 2m}\bn
    &\qquad = b_{p,j}(z_1,\ldots,z_n; q)
      \left(\prod_{l=1}^{m-p} 
      \D_j(\alpha + 2(m-p-l))\Phi\right)\Bigg|_{g, \alpha + 2(m - p)}.
  \end{align*}
  Now applying (\ref{D-shift}) several consecutive times, we obtain
  \begin{align*}
    &\left(\prod_{l=1}^m \D_j(\alpha + 2(m - l)) \right.\bn 
    &\left.\qquad + \sum_{p=1}^{m} b_{p,j}(z_1,\ldots,z_n; q)
      \prod_{l=1}^{m-p} \D_j(\alpha + 2(m-p-l))\right)\Phi|_{g,\alpha}
      = 0
  \end{align*}
  which concludes the proof.
\epfv

\noindent {\small \sc Department of Mathematics, Rutgers University,
110 Frelinghuysen Rd., Piscataway, NJ 08854-8019}
\vspace{1em}

\noindent {\em E-mail address}: francesco.fiordalisi@rutgers.edu 


\begin{thebibliography}{FGST2}

\bibitem[AF]{AF} F.~ Anderson and K.~Fuller, {\emph Rings and Categories
of Modules}, Graduate Texts in Mathematics 13, Springer-Verlag,
New York, Heidelberg, Berlin (1991).

\bibitem[AM1]{AM1} D.~Adamovi\'c and A.~Milas, Logarithmic intertwining
operators and $\mathcal W(2, 2p-1)$-algebras, {\em Journal of Math.
Physics} {\bf 48}, 073503 (2007).

\bibitem[AM2]{AM2} D.~Adamovi\'c and A.~Milas, On the triplet vertex
  algebra $\mathcal W(p)$, {\em Adv. Math.} {\bf 217}, 2664-2699 (2008).

\bibitem[AM3]{AM3} D.~Adamovi\'c and A.~Milas, On $\mathcal W$-algebras
associated to $(2,p)$ minimal models and their representations,
{\em IMRN} {\bf 20}, 3896-3934 (2010).

\bibitem[AM4]{AM4} D.~Adamovi\'c and A.~Milas, An analogue of modular BPZ
equation in logarithmic (super)conformal field theory,
{\em Contemporary Mathematics} {\bf 497}, (2009) 1-17. 

\bibitem[Ar]{Ar} Y.~Arike, Some remarks on symmetric linear functions
and pseudotrace maps, {\em Proc. Japan Acad. Ser. A Math. Sci.} {\bf 86}
(2010), 119-124.

\bibitem[AN]{AN} Y.~Arike and K.~Nagatomo, Some remarks on pseudo-trace 
functions for orbifold models associated with symplectic fermions.
{\em Int. J. Math.} {\bf 24}, (2013) 1350008 

\bibitem[B]{B} R.~E.~Borcherds, Vertex algebras, Kac-Moody algebras, and the
  Monster, {\em Proc. Natl. Acad. Sci.} {\bf 83}, 3068-3071 (1986).

\bibitem[Br]{Br} M.~Brou\'e, Higman criterion revisited, {\em Mich. J.
Math.} {\bf 58}, (2009) 125-179.

\bibitem[DLM1]{DLM1} Y.~Dong, H.~Li and G.~Mason, Vertex operator
algebras and associative algebras, {\em J. Algebra} {\bf 206}, 67-96
(1998).

\bibitem[DLM2]{DLM2} Y.~Dong, H.~Li and G.~Mason, Modular invariance
of trace functions in orbifold theory and generalized moonshine,
{\em Comm. Math. Phys.} {\bf 214}, 1-56 (2000).

\bibitem[FH]{FH} F. Fiordalisi and Y.-Z. Huang, Modular invariance for
logarithmic intertwining operators, in preparation.

\bibitem[FHL]{FHL} I.~Frenkel, Y.-Z.~Huang and J.~Lepowsky, On axiomatic
approaches to vertex operator algebras and modules, {\em Memoirs
Amer. Math. Soc.}, {\bf 104,} 1993.

\bibitem[FLM]{FLM} I.~Frenkel, J.~Lepowsky and A.~Meurman, {\em Vertex
Operator Algebras and the Monster}, Pure and Appl. Math., Vol. 134,
Academic Press, New York, 1988.

\bibitem[H1]{H1} Y.-Z.~Huang, {\em Two-dimensional conformal geometry and
vertex operator algebras}, Progress in Mathematics, Vol. 148,
Birkh\"auser, Boston, 1997.

\bibitem[H2]{H2} Y.-Z.~Huang, Differential equations, duality and modular
invariance, {\em Comm. Contemp. Math.} {\bf 7} (2005), 649-706.

\bibitem[H3]{H3} Y.-Z.~Huang, A theory of tensor product
for module categories for a vertex operator algebra, IV, {\em
J. Pure Appl. Algebra} {\bf 100} (1995), 173-216.

\bibitem[H4]{H4} Y.-Z.~Huang, Vertex operator algebras and the Verlinde
conjecture, {\em Comm. Contemp. Math.} {\bf 10} (2008), 103--154.

\bibitem[H5]{H5} Y.-Z.~Huang, Rigidity and modularity of vertex tensor
  categories, {\em Comm. Contemp. Math.} {\bf 10} (2008), 871--911.

\bibitem[H6]{H6} Y.-Z. Huang, Cofiniteness conditions, projective covers and 
the logarithmic tensor product theory, {\it J. Pure Appl. Alg.} {\bf 213} (2009), 458--475. 

\bibitem[HL1]{HL1} Y.-Z.~Huang and J.~Lepowsky, A theory of tensor product
for module categories for a vertex operator algebra, I, {\em
Selecta Mathematica (New Series)} {\bf 1} (1995), 699-756.

\bibitem[HL2]{HL2} Y.-Z.~Huang and J.~Lepowsky, A theory of tensor product
for module categories for a vertex operator algebra, II, {\em
Selecta Mathematica (New Series)} {\bf 1} (1995), 757-786.

\bibitem[HL3]{HL3} Y.-Z.~Huang and J.~Lepowsky, A theory of tensor product
for module categories for a vertex operator algebra, III, {\em
J. Pure Appl. Algebra} {\bf 100} (1995), 141-171.

\bibitem[HLZ1]{HLZ1} Y.-Z.~Huang, J.~Lepowsky and L.~Zhang, Logarithmic
tensor category theory for generalized modules for a conformal
vertex algebra, I: Introduction and strongly graded
algebras and their generalized modules, in:
{\em Conformal Field Theories and Tensor Categories, Proceedings of a Workshop
Held at Beijing International Center for Mathematics Research},
ed. C. Bai, J. Fuchs, Y.-Z. Huang, L. Kong, I. Runkel and C. Schweigert,
Mathematical Lectures from Beijing University, Vol. 2, Springer, New York, 2014,
169-248.  

\bibitem[HLZ2]{HLZ2} Y.-Z.~Huang, J.~Lepowsky and L.~Zhang, Logarithmic
tensor category theory, II: Logarithmic formal calculus and properties of 
logarithmic intertwining operators, arXiv:1012.4196.

\bibitem[HLZ3]{HLZ3} Y.-Z.~Huang, J.~Lepowsky and L.~Zhang, Logarithmic
tensor category theory, III: Intertwining maps and tensor product bifunctors,
arXiv:1012.4197.

\bibitem[HLZ4]{HLZ4} Y.-Z.~Huang, J.~Lepowsky and L.~Zhang, Logarithmic
tensor category theory, IV: Constructions of tensor product bifunctors and
the compatibility conditions, arXiv:1012.4198.

\bibitem[HLZ5]{HLZ5} Y.-Z.~Huang, J.~Lepowsky and L.~Zhang, Logarithmic
tensor category theory, V: Convergence condition for intertwining maps and the
corresponding compatibility conditions, arXiv:1012.4199.

\bibitem[HLZ6]{HLZ6} Y.-Z.~Huang, J.~Lepowsky and L.~Zhang, Logarithmic
tensor category theory, VI: Expansion condition, associativity of
logarithmic intertwining operators, and the
associativity isomorphisms, arXiv:1012.4202.

\bibitem[HLZ7]{HLZ7} Y.-Z.~Huang, J.~Lepowsky and L.~Zhang, Logarithmic
tensor category theory, VII: Convergence and extension properties and
applications to expansion for intertwining
maps, arXiv:1110.1929.

\bibitem[HLZ8]{HLZ8} Y.-Z.~Huang, J.~Lepowsky and L.~Zhang, Logarithmic
tensor category theory, VIII: Braided tensor category structure on
categories of generalized modules for a
conformal vertex algebra, arXiv:1110.1931.

\bibitem[HY]{HY} Y.-Z.~Huang and J.~Yang, Logarithmic intertwining
operators and associative algebras, {\em J. Pure Appl. Alg.} {\bf 216}
(2011), 1467-1492.

\bibitem[L]{L} S.~Lang, {\em Elliptic functions}, Graduate Texts in
Mathematics, Vol. 112, Springer-Verlag, New York, 1987.

\bibitem[LL]{LL} J.~Lepowsky and H.~Li, {\em Introduction to Vertex
Operator Algebras and Their Representations}, Progress in Mathematics,
VOl. 227, Birkh\"auser, Boston, 2003.

\bibitem[M1]{M1} A.~Milas, Weak modules and logarithmic intertwining operators
for vertex operator algebras, in {\em Recent Developments in
Infinite-Dimensional Lie Algebras and Conformal Field Theory},
ed. S. Berman, P. Fendley, Y.-Z. Huang, K. Misra, and B.
Parshall, Contemp. Math., Vol. 297, American Mathematical Society, Providence,
RI, 2002, 201–225. 

\bibitem[M2]{M2} A.~Milas, Logaritmic intertwining operators and vertex
operators, {\em Comm. Math.Phys.} {\bf 277} (2008), 497-529.

\bibitem[Miy1]{Miy1} M.~Miyamoto, Intertwining operators and modular
invariance, preprint, arXiv:math.QA/0010180

\bibitem[Miy2]{Miy2} M.~Myiamoto, Modular invariance of vertex operator
algebras satisfying $C_2$-cofiniteness, {\em Duke Math. J.} {\bf 122}
(2004), 51-91.

\bibitem[MS1]{MS1} G.~Moore and N.~Seiberg, Classical and quantum conformal
field theory, {\em Comm. Math. Phys.} {\bf 123} (1989), 177-254.

\bibitem[MS2]{MS2} G.~Moore and N.~Seiberg, Polynomial equations for rational
conformal field theories, {\em Phys. Lett.} {\bf B212} (1988), 451-460.

\bibitem[W]{W} C. ~Weibel {\em An Introduction to Homological Algebra},
Cambridge Studies in Adv. Math., Vol 38, Cambridge University Press,
Cambridge, 1994.

\bibitem[Z]{Z} Y.~Zhu, Modular invariance of characters of vertex
operator algebras, {\em J. Amer. Math. Soc.} {\bf 9} (1996), 237-307.


\end{thebibliography}
\end{document}